\documentclass[11pt]{amsart}

\usepackage{amsmath,amsfonts,amssymb,amsthm,enumerate,tikz-cd}

\usepackage{hyperref} 
\hypersetup{
	colorlinks=true, 
	     urlcolor= blue,  
	     linkcolor= blue, 
     citecolor=blue,	  
	bookmarksopen=true,           
}
\include{macros}

\def\Z{\mathbb Z}
\def\R{\mathbb R}
\def\Q{\mathbb Q}
\def\F{\mathbb F}

\def\Aut{\mathrm{Aut}}

\theoremstyle{plain}
\newtheorem{theorem}{Theorem}
\newtheorem{conjecture}{Conjecture}
\newtheorem{lemma}{Lemma}
\newtheorem{proposition}{Proposition}
\newtheorem{corollary}{Corollary}

\theoremstyle{definition}

\newtheorem{definition}{Definition}

\newtheorem{remark}{Remark}

\DeclareMathOperator{\ab}{ab}

\newcommand{\Coker}{\mathrm{Coker}}

\DeclareMathOperator{\Gal}{Gal}

\newcommand{\Tr}{\mathrm{Tr}}

\newcommand{\Roots}{\mathrm{Roots}}

\newcommand{\Cl}{\mathrm{Cl}}

\newcommand{\ur}{\mathrm{ur}}
\newcommand{\Pic}{\mathrm{Pic}}

\newcommand{\dR}{\mathrm{dR}}

\newcommand{\Ind}{\mathrm{Ind}}

\newcommand{\Sym}{\mathrm{Sym}}
\newcommand{\cris}{\mathrm{cris}}

\newcommand{\sep}{\mathrm{sep}}
\newcommand{\Sel}{\mathrm{Sel}}
\newcommand{\Nm}{\mathrm{Nm}}
\newcommand{\Jac}{\mathrm{Jac}}
\newcommand{\et}{\mathrm{\acute{e}t}}
\newcommand{\loc}{\mathrm{loc}}

\DeclareMathOperator{\rk}{rk}

\newcommand{\Res}{\mathrm{Res}}

\newcommand{\val}{\mathrm{val}}

\newcommand{\Ker}{\mathrm{Ker}}
\newcommand{\gr}{\mathrm{gr}}
\newcommand{\Hom}{\mathrm{Hom}}

\DeclareMathOperator{\ext}{Ext}

\begin{document}
\title{$2$-descent for Bloch--Kato Selmer groups and rational points on hyperelliptic curves I}
\author{Netan Dogra}
\maketitle
\pagestyle{headings}
\markright{$2$-DESCENT FOR BLOCH--KATO SELMER GROUPS I}

\begin{abstract}This paper introduces explicit Galois cohomological methods for determining the ranks of Bloch--Kato Selmer groups associated to the Tate twists of the 2-adic second \'etale cohomology of the Jacobian of a hyperelliptic curve with a rational Weierstrass point. In particular, this can give a method to determine the rational points on such curves via the Chabauty--Coleman--Kim method. This is applied to answer a question of Bugeaud, Mignotte, Siksek, Stoll and Tengely.
\end{abstract}

\tableofcontents

\section{Introduction}

Among curves of genus greater than one, hyperelliptic curves have traditionally been the most amenable to explicit determination of rational or integral points. For \textit{integral} points, i.e. for integer solutions to a given equation of the form
\[
y^2 =f(x),
\]
the effective computability was established by Baker \cite{baker}. For rational points, effective computability is still an open problem, but the Chabauty--Coleman method is often practical for determining the set of rational points. However both of these methods have drawbacks which limit their applicability for practical computation. For Baker's method, the limitation is that of determining integer points below the height bound. For the Chabauty--Coleman method, one of the main limitations is the condition on the Mordell--Weil rank.

For example, in \cite{BMSST}, Bugeaud, Mignotte, Siksek, Stoll and Tengely consider the genus two curve
\[
X: y^2 -y =x^5 -x.
\]
They combine refinements of Baker's method with Mordell--Weil sieving to determine the set of integer solutions. In the same paper, they ask whether the rational points can be determined. The Mordell--Weil rank of the Jacobian of $X$ is three, and the endomorphism algebra is $\mathbb{Z}$, so the Chabauty--Coleman method does not apply, and nor does the (usual) quadratic Chabauty method.

The present paper is devoted to developing methods which allow the computational verification of finiteness of $X(\Q _2 )_2 $, and which should lead to a practical method for computing rational points on many hyperelliptic curves of Mordell--Weil rank outside the Chabauty range.
\begin{theorem}\label{thm:BMSST_rational}
\[
X(\Q )=\left\{ \begin{array}{c} \infty ,(0, 1), (\frac{1}{4}, \frac{15}{32}), (2, 6),(3, -15), (1, 1), (30, -4929),  \\ (-1, 1), (1, 0), (30, 4930), (3, 16 ), (\frac{1}{4}, \frac{17}{32}), (2, -5), \\ (0, 0), (-1, 0), (-\frac{15}{16}, -\frac{185}{1024}), (-\frac{15}{16}, \frac{1209}{1024}) \end{array} \right\}
\]
\end{theorem}

\subsection{$2$-descent for hyperelliptic curves}
Let $K$ be a field of characteristic different from 2 and let $X/K$ be a smooth hyperelliptic curve given by a Weierstrass equation of the form
\[
y^2 =f(x)=\sum _{i=0}^{2g+1}a_i x^i ,
\]
with $a_{2g+1}\neq 0$. Then it is well-known that there is an isomorphism
\[
H^1 (K,J[2]) \simeq \Ker (K_f ^\times \otimes \mathbb{F}_2 \stackrel{\Nm }{\longrightarrow }K^\times \otimes \mathbb{F}_2 ),
\] 
where $K_f :=K[x]/(f(x))$ and $J$ is the Jacobian of $X/K$. Furthermore, as recalled below there is an explict field theoretic description of the Kummer map
\[
J(K)\otimes \F _2 \to H^1 (K,J[2]).
\]
When $K$ is a number field, this isomorphism may be used to describe the Selmer group $\Sel _2 (J)$ as a subgroup of $H^1 (K,J[2])$ (see \cite{schaefer}, \cite{PS97}, \cite{stoll:implementing}). This isomorphism is the basis of existing algorithms for proving upper bounds on the rank of the Jacobian of $X/K$, and is hence an essential tool in the implementation of the Chabauty--Coleman method for computing rational points on hyperelliptic curves whose Jacobians have small Mordell--Weil rank. 

In this paper, we develop methods to bound the rank of a different $2$-Selmer group, and give applications to computing the rational points on hyperelliptic curves whose Jacobians have small Mordell--Weil rank. Rather than using this finite Galois cohomology group to approximate the Mordell--Weil group of an abelian variety, we use it to approximate a Bloch--Kato Selmer group in the sense of \cite{BK}. Conjecturally, the dimension of this Selmer group should equal the rank of an associated higher Chow group, however such an implication is essentially never known outside the range of abelian varieties and $K$-groups of number fields.

In contrast to the mysterious phenomena arising when considering ranks of abelian varieties, for Galois representations $M$ associated to pure motives of weight $<-1$, Bloch and Kato conjecture that the dimension of the Selmer group $H^1 _f (\Gal (\overline{\Q }|\Q ),M)$ is essentially a purely \textit{geometric} invariant. More precisely, they conjecture a formula for it solely in terms of the action of complex conjugation on $M$, and the Hodge filtration on its de Rham realisation. There has been a vast amount of working on verifying this conjecture in various cases, which we will not attempt to survey here. However, both of the main approaches (via Euler systems, and via modularity lifting) depend on knowing that the underlying representation has an automorphic origin.


In this article, we propose a naive approach more in the spirit of computing the rank of an abelian variety via $2$-descent. Namely, given a $p$-adic representation $V$ of $\Gal (\Q _S |\Q )$ with a Galois stable lattice $T$, one may bound the $\Q _p $-dimension of $H^1 (\Gal (\Q _S |\Q ),V)$ by the $\F _p $-dimension of $H^1 (\Gal (\Q _S |\Q ),T\otimes \F _p )$. Hence a first step is to be able to estimate the dimension of $H^1 (\Gal (\Q _S |\Q ),T\otimes \F _p )$. We give such an estimate when $p=2$ and $V=\wedge ^2 V_2 J$ where $J$ is the Jacobian of a hyperelliptic curve with a rational Weierstrass point, however the methods have a somewhat broader applicability. For example, if $J$ is as above then these methods give a (weaker) bound for the dimension of $H^1 (\Gal (\Q _S |\Q ),T\otimes \F _2 )$ when $T$ is \textit{any} Galois representation generated by $V_2 J$ (i.e. any summand of $V^{\otimes i}(j)$ for $i\geq 0$ and $j\in \mathbb{Z}$).

With more work, one can describe a subspace which one might call $H^1 _f (\Gal (\Q _S |\Q ),T\otimes \F _p )$ inside $H^1 (\Gal (\Q _S |\Q ),T\otimes \F _p )$ such that
\[
\dim H^1 _f (G_{\Q ,S},V)\leq \dim H^1 _ f (G_{\Q ,S},T\otimes \F _p )
\]
The notation $H^1 _f (\Gal (\Q _S |\Q ),T\otimes \F _p )$ is perhaps inadvisable, because in fact the subspace $H^1 _f (G_{\Q ,S},T\otimes \F _p )$ will \textit{not} be an invariant of $T\otimes \F _p $ but will depend on $T$. Recent work of Iovita and Marmora explores the extent to which the image of $H^1 _f (\Q _p ,T)$ in $H^1 (\Q _p ,T\otimes \F _p )$ depends on the lift $T$ of $T\otimes \F _p $ \cite{IM15}. In \cite{gazaki}, Gazaki explores this in the case where $T$ is a tensor product of Tate modules of abelian varieties, which is essentially the case of interest for quadratic Chabauty. In these papers one can obtain a subspace depending only on $T\otimes \F _p $ if one imposes a condition on the Hodge--Tate weights, but since we are interested in $p=2$ in this paper we are working outside the range where their results apply.

One tool used for studying the Bloch--Kato Selmer group $H^1 _f (\Gal (\overline{\Q }|\Q ),\wedge ^2 V)$ is its relation to quotients of the \'etale fundamental group. More precisely, let $K|\Q $ be a field extension, and let $b$ be a point in $X(K)$. Let $\pi _1 ^{\et ,2}(X_{\overline{K}}-\{ \infty \} ,b)$ denote the maximal pro-$2$ quotient of the \'etale fundamental group of $X_{\overline{K}}-\{ \infty \} $ with basepoint $b$. Let $\Pi _2 $ denote its maximal $2$-nilpotent quotient. Then $\Pi _2 $ is a central extension of $T_2 (J)$ by $\wedge ^2 T_2 (J)$.
There is a unique quotient $\Gamma $ of $\Pi _2 $ sitting in a commutative diagram with exact rows
\begin{equation}\label{eqn:gammabar}
\begin{tikzcd}
1 \arrow[r] & \wedge ^2 T_2 (J) \arrow[d] \arrow[r] & \Pi _2 \arrow[d] \arrow[r] & T_2 (J) \arrow[d] \arrow[r] & 1 \\
1 \arrow[r] & \wedge ^2 J[2] \arrow[r]           & \Gamma  \arrow[r]           & J[4] \arrow[r]           & 1
\end{tikzcd}
\end{equation}
where the outer vertical maps are the obvious ones. The main results on ``$\Gamma$-descent" for $X$ proved in this paper may be viewed as nonabelian generalisations of the main results described above for $2$-descent on hyperelliptic curves. Firstly, we give an explicit field-theoretic description of the cohomology group $H^1 (K,\wedge ^2 J[2])$, in terms of units mod squares of a certain \'etale algebra of degree $2g^2 +g$ over $K$. If $\Gal (\overline{K}|K)$ acts $2$-transitively on the roots of the polynomial $f$, then this \'etale algebra is a field extension of $K$.

Given a pair of points $z,b$ in $X(K)$, we obtain a class $\kappa _{\Gamma  ,b}(z)$ in $H^1 (K,\Gamma )$ whose image in $H^1 (K,J[4])$ is the mod $4$ Kummer class $\kappa _{4,b}(z)$ of $[z]-[b]$. If $[z]-[b]$ lies in $4\cdot \Jac (X)(K)$, there is an element of $H^1 (K,\wedge ^2 J[2])$ with image $\kappa _{\Gamma  ,b}(z)$ in $H^1 (K,\Gamma  )$. We give an explicit formula for such a class in the case when $X$ is hyperelliptic with a rational Weierstrass point, and the points $z$ and $b$ are not Weierstrass points.

To state the result, let $X/K$ be a hyperelliptic curve given by a Weierstrass equation of the form
\[
y^2 =f(x):=\sum _{i=0}^{2g+1}a_i x^i ,
\]
with $a_{2g+1}\neq 0$. 
Let $K_{f,2}$ denote the \'etale algebra
\[
K[s,t,\frac{1}{s-t}]/(f(s),f(t)),
\] 
and define $K_f ^{(2)}:=K_{f,2}^{\tau }$, where $\tau $ is the involution sending $s$ to $t$. Let $\alpha $ and $\beta $ denote the image of $s$ and $t$ in $K_{f,2}$.

\begin{theorem}
\begin{enumerate}
\item We have an isomorphism
\[
H^1 (K,\wedge ^2 J[2]) \simeq \Ker (K_f ^{(2),\times } \otimes \mathbb{F}_2 \stackrel{\Nm }{\longrightarrow }K_f ^\times \otimes \mathbb{F}_2 ).
\]
\item Let $b$ and $z$ be non-Weierstrass points of $X(K)$. Suppose that $\kappa _{b,4}(z)=0$. Then, with respect to the isomorphism in part (1) of the theorem, a lift of $\kappa _{b,\Gamma }(z)$ to $H^1 (K,\wedge ^2 J[2])$ is given by 
\[
2\left( 1+\frac{c_{\alpha }u_{\alpha }-c_{\beta }u_{\beta }}{\beta -\alpha } \right)
\]
where $c_{\alpha }=x(b)-\alpha $ and $u_{\alpha }$ is a square root of $(x(z)-\alpha )/(x(b)-\alpha )$.
\end{enumerate}
\end{theorem}
In the sequel to this paper \cite{BKdescent2}, we study liftability or `second descent' obstructions for $H^1 _f (\Q _v ,\wedge ^2 T_2 J)$, and use these to verify finiteness of $X(\Q _2 )_2 $ for several thousand genus 2 curves with Jacobians of Mordell--Weil rank 2 and Picard number 1. Recent work of Berry \cite{berry} has shown how to extend many of these results to the case of hyperelliptic curves of higher genus, and of even degree (i.e. hyperelliptic curves without a rational Weierstrass point). Berry's work also shows how the crystalline conditions, studied in section \ref{sec:2} using a nonabelian Kummer map, may be studied more systematically in the case where the curve has good ordinary reduction at 2.
\subsection{Notation}
For a number field $K$ and a finite set $S$ of places of $K$, we define $\Cl (\mathcal{O}_{K,S}):=\Pic (\mathcal{O}_{K,S})$. For any integer $N>0$, we define $(K^\times \otimes \mathbb{Z}/N\mathbb{Z})_S $ to be the subgroup of $K^\times \otimes\mathbb{Z}/N\mathbb{Z}$ consisting of elements whose image in $K_v ^\times \otimes \mathbb{Z}/N\mathbb{Z}$ lies in the image of $\mathcal{O}_v \otimes \mathbb{Z}/N\Z $ for all $v$ not in $S$.

When we say an algebraic group is $n$-unipotent, we mean that it is unipotent, and all $n$-fold commutators vanish (with the convention that a $1$-fold commutator is just a commutator). For example, a $1$-unipotent group is a unipotent group which is abelian.

Given an abelian group $M$, we define
\[
\wedge ^2 M:= M^{\otimes 2}/\langle x\otimes x |x\in M \rangle .
\]

Given a field $K$ and a $\Gal (K^{\sep }|K)$-module $M$, we will sometimes write $H^i (K,M)$ for the Galois cohomology group $H^i (\Gal (K^{\sep }|K),M)$.
\subsection*{Acknowledgements}
I would like to thank the anonymous referees for helpful corrections and suggestions for changes to the structure of the paper. I am grateful to Ed Schaefer for explaining $2$-Selmer groups of hyperelliptic curves to me, and to Michael Stoll for several detailed suggestions about an earlier version of this work, and corrections and comments on an earlier version of this paper. I am grateful to Lee Berry for many corrections to a previous version of this paper and useful discussions. This research was supported by a Royal Society University Research Fellowship.
\section{Galois cohomology of the wedge square}\label{sec:shapiro}
In this section, unless otherwise indicated, $K$ will denote a field of characteristic different from $2$ and $f\in K[x]$ will be a monic separable polynomial of degree $d$. None of our general results on \'etale algebras and Galois cohomology are original (see for example \cite{BPS} or \cite{SS}), but we include them for the sake of completeness and to establish notation.

First we define an essential inverse to the usual functor from finite \'etale $K$-algebras to finite $\Gal (K)$-sets. 
Given a finite $\Gal (K^{\sep }|K)$-set $S$. we define $K_S$ to be the \'etale algebra of $\Gal (K^{\sep }|K)$-equivariant maps (of sets) $S\to K^{\sep }$. If $S=\Gal (K)/H$ where $H<\Gal (K)$ is a finite index subgroup, then we have an isomorphism $K_S \simeq (K^{\sep })^H$ given by sending an equivariant map $f$ to $f(H)$.

Recall from the introduction the \'etale algebras $K_f ,K_{f,2}$ and $K_f ^{(2)}$. For an \'etale $K$-algebra $L$, with factorisation into fields $\prod _{i=1}^n L_i $, and a finite abelian group with an action of $\Gal (L_i )$ for each $i$, we define $\Ind ^L _K M:=\oplus _{i=1}^n \Ind ^{L_i }_K M$.

\begin{lemma}
Let $G$ be a group, let $S$ be a finite $G$-set, and let $k$ be a field of characteristic $2$. Then we have an isomorphism of $G$-modules
\[
\wedge ^2 k[S]\simeq k[S^{(2)}]
\]
where $S^{(2)}$ denotes the set of pairs of unordered distinct elements of $S$.
\end{lemma}
\begin{proof}
This is given by sending $[x]\wedge [y]$ to $[\{ x,y\} ]$.
\end{proof}

\begin{lemma}
Let $S$ denote the $\Gal (K)$-set $\Hom _K (K_f ,K^{\sep })$. Then we have an isomorphism
\[
K_{S^{(2)}}\simeq K_f ^{(2)}.
\]
\end{lemma}
\begin{proof}
It is enough to prove that we have an isomorphism of $\Gal (K)$-sets
\[
S^{(2)}\simeq \Hom _K (K_f ^{(2)},K^{\sep }).
\]
Let $S^2 -\Delta $ denote the set of ordered pairs of distinct elements of $S$. Then we have an isomorphism of $\Gal (K)$-sets
\[
S^2 -\Delta \simeq \Hom _K (K_{f,2},K^{\sep })
\]
given by sending $(\alpha ,\beta )$ to $\alpha \otimes \beta $. This isomorphism is equivariant with respect to the natural `swapping' involutions $\tau $ on both sides. Since we have an isomorphism of $G$-sets 
\[
(S^2 -\Delta )/\tau \simeq S^{(2)},
\]
it is enough to prove that pulling back by $K_f ^{(2)}\to K_{f,2}$ induces an isomorphism of $G$-sets
\begin{equation}\label{eqn:need_iso}
\Hom _K (K_{f,2},K^{\sep })/\tau \simeq \Hom _K (K_f ^{(2)},K^{\sep })
\end{equation}
The pullback
\[
\Hom _K (K_{f,2},K^{\sep })\to \Hom (K_f ^{(2)},K^{\sep })
\]
factors through quotienting by $\tau $, and is surjective. The isomorphism then follows from the fact that the orders of the right and left hand sides of \eqref{eqn:need_iso} are the same.
\end{proof}
Let $\Nm :\F _2 [S^{(2)}]\to \F _2 [S]$ denote the map sending $[\{\alpha ,\beta \}] $ to $[\alpha ] +[\beta ]$. Let
\[
\Nm :K_f ^{(2),\times }\otimes \F _2 \to K_f \otimes \F _2 
\]
denote the composite map
\[
K_f ^{(2),\times }\otimes \F _2 \to K_{f,2}^\times \otimes \F _2 \to K_f \otimes \F _2
\]
where the first map is induced by the inclusion $K_f ^{(2)}\to K_{f,2}$ and the second is the norm map associated to the map of \'etale algebras $K_f \to K_{f,2}$ sending $\overline{x}$ to $\overline{x}$.

\begin{lemma}
Given a field $K$ of characteristic different from $2$, and a finite $\Gal (K)$-set $S$, we have an isomorphism
\[
H^1 (K,\F _2 [S])\simeq K_S ^\times \otimes \F _2 
\]
\end{lemma}
\begin{proof}
It is enough to prove this when $S$ is a transitive $G$-set. In that case $S$ is isomorphic to $\Gal (K) \backslash \Gal (L)$ for a finite separable extension $L|K$, and $K_S $ is isomorphic to $L$ as above. On the other hand $\F_ 2[\Gal (K)\backslash \Gal (L)]$ is isomorphic to $\Ind ^L _K \F _2 $, and hence the statement follows from Shapiro's lemma.
\end{proof}

This identification satisfies some obvious functorial properties.
\begin{lemma}\label{lemma:the_obvious_one}
Let $\pi :S_1 \to S_2 $ be a surjective map of $G$-sets, and let $\pi ^*$ denote the induced map $\F _2 [S_2 ]\to \F _2 [S_1 ]$ sending $s_2 $ to $\sum _{\pi (s_1 )=s_2 }s_1 $. Then the map on cohomology
\[
K_{S_2 }^\times \otimes \F _2 \to K_{S_1 }^\times \otimes \F _2 
\]
is obtained from the natural map
\[
K_{S_2 }\to K_{S_1 }
\]
induced by composition of maps $S_1 \to S_2 \to K^{\sep }$.
\end{lemma}
\begin{proof}
We may reduce to the case where $G$ acts transitively on $S_1 $ and $S_2 $, and $S_i $ is isomorphic to $\Gal (K)/\Gal (L_i )$ for $L_i |K$ finite extensions of $K$. Then $\pi $ corresponds to a unique map of $K$-algebras $L_2 \to L_1 $, and $\pi $ is isomorphic to the induction from $\Gal (L_2 )$ to $\Gal (K)$ of the map of $\Gal (L_2 )$-sets
\begin{equation}\label{eqn:induced_res}
\F_2 \to \F _2 [\Gal (L_2 )/\Gal (L_1 )]
\end{equation}
sending $1$ to $\sum _{x\in \Gal (L_2 )/\Gal (L_1 )}x$. The Lemma now follows from Shapiro's lemma, since \eqref{eqn:induced_res} induces the natural map $L_2 ^\times \otimes \F _2 \to L_1 ^\times \otimes \F _2 $ on cohomology.
\end{proof}

\begin{lemma}\label{lemma:norm_is_norm}
Let $S_1 $ and $S_2 $ be $G$-sets. Let 
\[
\Nm :\F _2 [S_1 \times S_2 ]\to \F _2 [S_2 ]
\]
be the norm map.
We have a commutative diagram
\[
\begin{tikzcd}
H^1 (K ,\F _2 [S_1 \times S_2 ]) \arrow[r] \arrow[d] & H^1 (K ,\F _2 [S_2 ]) \arrow[d] \\
(K_{S_1 }\otimes K_{S_2 })^\times \otimes \F _2 \arrow[r] &  K_{S_2 }^\times \otimes \F _2 .          
\end{tikzcd}
\]
\end{lemma}
\begin{proof}
We may reduce to the case where $G$ acts transitively on $S_2 $. Suppose $s\in S_2 $ is stabilised by $\Gal (L)$. Then we want to identify the norm map with the map
\[
L_{S_1 }^\times \otimes \F _2 \to L^\times \otimes \F _2 .
\]
The norm map $\F _2 [S_1 \times S_2 ]\to \F _2 [S_2 ]$ is the induction to $\Gal (K)$ of the usual norm map of $\Gal (L)$ modules
\[
\F _2 [S_1 ]\to \F _2 ,
\]
and hence the result follows from Shapiro's lemma.
\end{proof}

\begin{lemma}\label{lemma:random}
We have a commutative diagram
\[
\begin{tikzcd}
H^1 (K ,\F _2 [S^{(2)}]) \arrow[r] \arrow[d] & H^1 (K_2 ,\F _2 [S]) \arrow[d] \\
K_f ^{(2),\times } \otimes \F _2  \arrow[r]           & K_{f}^\times \otimes \F _2 .          
\end{tikzcd}
\]
\end{lemma}
\begin{proof}
The norm map $\Nm$ on Galois modules is the composite 
\begin{equation}\label{eqn:norm_is_composite}
\F _2 [S^{(2)}]\to \F _2 [S ]\otimes \F _2 [S] \to \F _2 [S]
\end{equation}
where the first map sends $\{\alpha ,\beta \}$ to $\alpha \otimes \beta +\beta \otimes \alpha $, and the send map sends $\alpha \otimes \beta $ to $\beta $. Hence the lemma follows from commutativity of 
\[
\begin{tikzcd}
H^1 (K ,\F _2 [S^{(2)}]) \arrow[r] \arrow[d] & H^1 (K,\F _2 [S\times S]) \arrow[r] \arrow[d] & H^1 (K_2 ,\F _2 [S]) \arrow[d] \\
K_f ^{(2),\times } \otimes \F _2  \arrow[r] & (K_f \otimes K_f )^\times \otimes \F _2 \arrow[r] & K_{f}^\times \otimes \F _2 .          
\end{tikzcd}
\]
where the right hand square is as in Lemma \ref{lemma:norm_is_norm}, the top line is obtained from \eqref{eqn:norm_is_composite}, and the bottom left horizontal arrow is induced from the inclusion of $K_f ^{(2)}$ into $K_{f,2}$ and the isomorphism $K_{f,2}\times K_f \simeq K_f \otimes K_f $. Hence we reduce to commutativity of the left-hand square, which is a special case of Lemma \ref{lemma:the_obvious_one}.

\end{proof}

We henceforth suppose $f$ has odd degree, and define a sub-object $J[2]$ of $\Ind ^{K_f }_K \F _2 $ as $J[2]:=\Ker (\Ind ^{K_f }_K \F _2 \stackrel{\Nm}{\longrightarrow } \F _2 )$ (the notation is justified by the fact that $J[2]$ is indeed isomorphic to the $2$-torsion of the Jacobian of the hyperelliptic curve defined by $f$). The short exact sequence
\begin{equation}\label{eqn:seq}
0\to J[2]\to \Ind ^{K_f }_K \F _2 \to \F _2  \to 0
\end{equation}
has a Galois equivariant splitting, given by sending $1\in \F _2 $ to $\sum _{\alpha \in \Roots (f)}\alpha $. We deduce the following proposition.

\begin{proposition}\label{prop:main1}
Let $f$ have odd degree. Then we have an isomorphism
\[
H^1 (K,\wedge ^2 J[2])\simeq \Ker ((K_f ^{(2)})^\times \otimes \F _2 \to K_f ^\times \otimes \F _2 ).
\]
The image of the norm map $(K_f ^{(2)})^\times \otimes \F _2 \to K_f ^\times \otimes \F _2 $ is equal to $\Ker (K_f ^\times \otimes \F _2 \stackrel{\Nm }{\longrightarrow} K^\times \otimes \F _2 )$. 
\end{proposition}
\begin{proof}
The short exact sequence \eqref{eqn:seq} induces a short exact sequence
\begin{equation}\label{eqn:seq2}
0\to \wedge ^2 J[2]\to \wedge ^2 \Ind ^{K_f }_K \F _2 \to J[2]\to 0.
\end{equation}
The splitting of \eqref{eqn:seq} induces a splitting of \eqref{eqn:seq2}, inducing an isomorphism
\[
H^1 (K,\wedge ^2 J[2]) \simeq \Ker (H^1 (K,\wedge ^2 \Ind ^{K_f }_K \F _2 )\to H^1 (K,J[2]))
\]
On the other hand the splitting of \eqref{eqn:seq} tells us that this is isomorphic to $\Ker (H^1 (K,\wedge ^2 \Ind ^{K_f }_K \F _2 )\to H^1 (K,\Ind ^{K_f }_K \F _2 ))$. Here the map is equal to the map on cohomology induced by 
\[
\F _2 [\Roots (f)^{(2)}]\to \F _2 [\Roots (f)]
\]
from Lemma \ref{lemma:random}. Hence Lemma \ref{lemma:random} gives an identification
\[
H^1 (K,\wedge ^2 J[2])\simeq \Ker (K_f ^{(2),\times }\otimes \F _2 \to K_f ^\times \otimes \F _2 )
\]
with the kernel of the norm map. The splitting of \eqref{eqn:seq2} implies that $H^1 (K,\wedge ^2 \Ind ^{K_f }_K \F _2 )$ surjects onto $H^1 (K,J[2])$, and hence we deduce that the norm map surjects onto $\Ker (K_f ^{\times }\otimes \F _2 \to K^\times \otimes \F _2 )$.
\end{proof}


\section{The nonabelian $(x-T)$-map}
In this section we give what may be viewed as a `nonabelian generalisation' of the $(x-T)$ map defined by Schaefer \cite{schaefer} generalising a construction of Cassels \cite{cassels}.

First, for any field $K$, geometrically connected $K$-scheme $Z$, and $K$-points $x,y,$ we define $j_{Z,x}(y)\in H^1 (K,\pi _1 ^{\et }(Z_{K^{\sep }},x))$ to be the cohomology class of $\pi _1 ^{\et }(Z_{K^{\sep }};x,y)$.

Now let $K$ and $f$ be as in section \ref{sec:shapiro}, and suppose $f$ has odd degree. Let $X$ denote the hyperelliptic curve associated to $f$, and suppose $b,z\in X(K)$ are non-Weierstrass points. 
We have a nonabelian Kummer map
\[
(X-\{\infty \})(K) \stackrel{j_{X-\{ \infty \} ,b}}{\longrightarrow }H^1 (K,\pi _1 ^{\et }(X_{K^{\sep }}-\{\infty \} ,b)) 
\]
which may be defined by sending a point $z\in X-\{ \infty \} (K)$ to the Galois equivariant continuous $\pi _1 ^{\et}(X_{K^{\sep }}-\{ \infty \} ,b)$-torsor of paths $\pi _1 ^{\et }(X_{K^{\sep }}-\{ \infty \} ;b,z)$ (and similarly for $X-\{\infty \}$). More generally, for any geometrically connected scheme $Z$ over $K$, and $b\in Z(K)$, and any $\Gal (K^{\sep }|K)$-stable quotient $G$ of $\pi _1 ^{\et }(Z_{K^{\sep }},b)$, we define a map
\[
j_{G,b}:Z(K)\to H^1 (K,G)
\]
to be the composite of $j_{Z,b}$ with the pushforward map
\[
H^1 (K,\pi _1 ^{\et }(Z_{K^{\sep }},b))\to H^1 (K,G).
\] 

For example, let $\Gamma$ denote the quotient defined in \eqref{eqn:gammabar}, and let $\kappa _{\Gamma  ,b}(z)$ in $H^1 (K,\Gamma)$ denote the section class associated to $b$ and $z$. Then the image of $\kappa _{\Gamma  ,b}(z)$ in  $H^1 (K,J[4])$ is the mod $4$ Kummer class $\kappa _{4,b}(z)$ of $[z]-[b]$. If $[z]-[b]$ lies in $4\cdot \Jac (X)(K)$, there is an element of $H^1 (K,\wedge ^2 J[2])$ with image $\kappa _{\Gamma  ,b}(z)$ in $H^1 (K,\Gamma  )$.

For any root $\gamma $ of $f$ over an \'etale $K$-algebra $L$, we let $c_{\gamma }:=x(b)-\gamma $. If $z\in X(K)$ is a non--Weierstrass point with the property that $z-b$ lies in $2\cdot J(K)$, then there is an element $u_{\gamma }$ of $L$ such that $c_{\gamma }u_{\gamma }^2 =x(z)-\gamma $. For example, we can apply this construction to the \'etale algebra $K_{f,2}$, to get elements $c_{\alpha },u_{\alpha },c_{\beta }$ and $u_{\beta }$. The element $\frac{c_{\alpha }u_{\alpha }-c_{\beta }u_{\beta }}{\beta -\alpha }$ will in fact be defined over the field $K_f ^{(2)}$. This element of $K_f ^{(2)}$ is related to $j_{J[4],b}(z)$ and $j_{\Gamma ,b}(z)$ by the following `nonabelian $(x-T)$ maps'.
\begin{proposition}\label{prop:main2}
\begin{enumerate}
\item Suppose $z\in X(K)$ is a non-Weierstrass point with the property that $j_{J[2],b}(z)=0$ in $H^1 (K,J[2])$. Then a lift of $j_{J[4] ,b}(z)$ to $H^1 (K,J[2])$, via the exact sequence
\[
H^1 (K,J[2])\to H^1 (K,J[4])\to H^1 (K,J[2]),
\] 
is given by 
\[
\Nm _{K_f ^{(2)}|K_f }\left(1+ \frac{c_{\alpha }u_{\alpha }-c_{\beta }u_{\beta }}{\beta -\alpha }\right) .
\]
\item Suppose $z\in X(K)$ is a non-Weierstrass point with the property that $j_{J[4],b}(z)=0$ in $H^1 (K,J[4])$. Then a lift of $j_{\Gamma  ,b}(z)$ to $H^1 (K,\wedge ^2 J[2])$ is given by 
\[
2\left(1+ \frac{c_{\alpha }u_{\alpha }-c_{\beta }u_{\beta }}{\beta -\alpha }\right) .
\]
\end{enumerate}
\end{proposition}
The proof of this proposition will occupy the remainder of this section. We first translate this into the problem of constructing explicit $K$-models for certain covers of $X$, and describing their fibres above the point $z$. Let $X/K$ be a geometrically connected scheme over a field $K$. Let $f:Y\to X_{\overline{K}}$ be a finite Galois cover with group $G$. Then we have a surjection
\begin{equation}\label{eqn:pi1_surj}
\pi _1 ^{\et }(X_{\overline{K}} ,b)\to G.
\end{equation}
The following lemma is a special case of \cite[Corollary 34]{stix:arithmetic}.
\begin{lemma}\label{lemma:descent}
Suppose the kernel of $\pi _1 ^{\et }(X_{\overline{K}},b)\to G$ is stable under the conjugation action of $\Gal (\overline{K}|K)$.
\begin{enumerate}
\item There is a descent $(Y',c)$ of $Y$ to a pointed cover $(Y',c)\to (X,b)$ defined over $K$. 
\item The $G$-torsor structure on $f^{-1}(z)$ is $\Gal (\overline{K}|K)$-equivariant with respect to the Galois action on $G$ induced by \eqref{eqn:pi1_surj} and that on $f^{-1}(z)$ induced by the descent of $f$ to $K$.
\item The class of $f^{-1}(z)$ in $H^1 (K,G)$ is equal to the image of the class of $z$ in $H^1 (K,\pi _1 ^{\et }(X_{\overline{K}},b))$ in $H^1 (K,G)$ under the pushforward induced by \eqref{eqn:pi1_surj}.
\end{enumerate}
\end{lemma}

To apply this construction, we need a way of constructing $K$-models of the covers of $X_{K^{\sep }}$ corresponding to $J[2],J[4]$ and $\Gamma$. We will do this using some properties of the Weil restriction functor \cite[7.6]{BLR}. The construction we give below can also be found in \cite[\S 5]{stoll:descent} in a more general context. Recall that, given an \'etale algebra $L$ over a field $K$, and a scheme $X$ over $L$, a Weil restriction $\Res _{L|K}X$ is a $K$-scheme representing the functor on $K$-schemes
\[
T\mapsto \Hom _L (T\times _K L,X).
\]
In particular, we have an isomorphism $\Res _{L|K}(X)(K)\simeq X(L)$. The existence of a Weil restriction in the context in which we apply it will be a consequence of \cite[Theorem 7.6.4]{BLR}.

To descend $\mathbb{Z}/2\mathbb{Z}$ covers, we will make use of the following well-known calculation.
\begin{lemma}\label{lemma:sqrt_nm}
Let $K$ be a field of characteristic different from $2$, and let $K(\sqrt{d})|K$ be a quadratic \'etale algebra. Suppose 
\[
\Nm _{K(\sqrt{d})|K}(x+y\sqrt{d})=z^2 ,
\]
for $x,y,z\in K$, and $z$ nonzero. Then 
\[
x+y\sqrt{d}=\frac{(x+z+y\sqrt{d})^2}{2(x+z)}.
\]
In particular,
\[
(x+y\sqrt{d})\equiv 2(x+z)
\]
in $K(\sqrt{d})^\times \otimes \F _2 $.
\end{lemma}
Note that this formula is independent of the choice of a square root of the norm of $x+y\sqrt{d}$, since $x^2 -\Nm (x+y\sqrt{d})$ is a square in $K(\sqrt{d})$.


Let $f\in K[x]$ be a separable polynomial of odd degree and let $X$ be the associated hyperelliptic curve. Suppose we have a non-Weierstrass point $b\in X(K)$. Then we have a $\mathbb{Z}/2\mathbb{Z}$ cover of smooth curves over $K_f $
\begin{equation}\label{eqn:Yalpha}
Y_{\alpha }\to X_{K_f }
\end{equation}
given on function fields by $K_f (Y_{\alpha })=K_f (X)(\sqrt{c_{\alpha }(x-\alpha )})$. Recall that $c_{\alpha }:=x(b)-\alpha $, so in particular $b\in X(K_f )$ lifts to a $K_f $-point of $X$. Explicitly, $Y_{\alpha }$ is a hyperelliptic curve given by
\[
v_{\alpha }^2 =c\cdot c_{\alpha }\cdot \prod _{\beta \neq \alpha }(u_{\alpha }^2 -(\beta -\alpha )/c_{\alpha }),
\]
with map
\[
Y_{\alpha }\to X_{K_f}
\]
given by 
\[
(u_{\alpha } ,v_{\alpha })\mapsto (c_{\alpha }u_{\alpha }^2 +\alpha ,c_{\alpha }^g v_{\alpha }u_{\alpha }).
\]
The point $b$ lifts to the point $(1,y(b)/c_\alpha ^g )$.

Given a morphism $f:X\to Y$ of $L$-schemes admitting Weil restrictions to $K$, we obtain a morphism
\[
\Res _{L|K}(f):\Res _{L|K}X\to \Res _{L|K}Y
\]
defined as follows. Let 
\[
\theta :(\Res _{L|K}X )_L \to X
\]
be the morphism adjoint to the identity morphism on $X$. Then we define $\Res _{L|K}(f)$ to be the morphism adjoint to 
\[
f\circ \theta : (\Res _{L|K}X)_L \to Y.
\]
In particular, taking Weil restrictions from $K_f$ to $K$ we may descend the cover $Y_{\alpha }\to X_{K_f }$ above to a morphism of $K$-schemes
\[
\Res _{K_f |K}Y_{\alpha }\to \Res _{K_f |K}X_{K_f }.
\]
On the other hand we have a natural map $i:X\to \Res _{K_f |K}X_{K_f }$, and hence we may define $X_2 $ to be the fibre product $X\times _{\Res _{K_f |K}X_{K_f }}\Res _{K_f |K}Y_{\alpha }.$

\begin{lemma}\label{lemma:WeilRes}
Let $L|K$ be a finite extension, and let $M|L$ be an extension of the Galois closure of $L|K$. Let $X$ be an $L$-scheme admitting a Weil restriction to $K$. Then there is an isomorphism of $M$-schemes
\begin{equation}\label{eqn:WeilresGalois}
(\Res _{L|K}X)_M \simeq \prod _{\sigma :L\to M}X_M .
\end{equation}
where the product is over all $K$-embeddings of $L$ into $M$. Moreover, if $f:X\to Y$ is a morphism of $L$-schemes and $X$ and $Y$ admit a Weil restriction to $K$, then we have a commutative diagram
\[
\begin{tikzcd}
(\Res _{L|K}X)_M \arrow[d, "\simeq"] \arrow[r, "(\Res _{L|K}f)_M"] & (\Res _{L|K}Y)_M \arrow[d, "\simeq "] \\
\prod _{\sigma }X_M \arrow[r, "\prod _{\sigma :L\to M}f_{\sigma }"]                                                     & \prod _{\sigma }Y_M                                                
\end{tikzcd}
\]
\end{lemma}
\begin{proof}
This is a consequence of the corresponding identities between the functors which the respective Weil restrictions represent, see \cite[7.6, p.192]{BLR}.
\end{proof}

\begin{lemma}\label{lemma:X2}
The curve $X_2 $ is a product of two geometrically irreducible curves $X_2 ^+ $ and $X_2 ^-$, each of which is an \'etale $(\mathbb{Z}/2\mathbb{Z})^{2g}$ cover corresponding to $J[2]$.
\end{lemma}
\begin{proof}
We have a norm map
\[
\Res _{K_f |K}\mathbb{P}^1 \to \mathbb{P}^1
\]
and we obtain by composition with the map to $\Res _{K_f |K}X_{K_f }$ and the map $\Res (u_{\alpha }):\Res _{K_f |K}X_{K_f }\to \Res _{K_f |K}\mathbb{P}^1 $ a map
\[
\Nm (u_{\alpha }):X_2 \to \mathbb{P}^1 .
\]
We claim that $X_2 $ is a disjoint union of connected components $X_2 ^+ =\{ \Nm (u_{\alpha })\}=y$ and $X_2 ^- = \{ \Nm (u_{\alpha })=-y \}$. 
It is enough to prove this after base change to $K^{\sep }$. The function field of the covering corresponding to $J[2]$ is generated over $K^{\sep }(X)$ by $\{ u_{\alpha }:\alpha \in \Roots (f) \}$. The result then follows from base changing to $K^{\sep }$ and applying Lemma \ref{lemma:WeilRes}.
\end{proof}
By construction, $b\in X(K)$ lifts to a $K$-point $\widetilde{b}$ of $X_2 ^+$, and hence we obtain a pointed cover of $(X,b)$ defined over $K$.

\begin{lemma}\label{lemma:Z4Zcovers}
For each $\beta \neq \alpha $ a root of $f$, let $(\beta -\alpha )^{1/2}$ denote a choice of square root of $\beta -\alpha $ in $K^{\sep }$. Let
\[
g_{\alpha ,0}:=\prod _{\beta \neq \alpha }(u_{\alpha }-(\beta -\alpha )^{1/2})
\] 
Then $K^{\sep }(X)(u_{\alpha })(\sqrt{g_{\alpha ,0}})$ is a $\mathbb{Z}/4\mathbb{Z}$-extension of $K^{\sep }(X)$ extending to an \'etale cover of $X_{K^{\sep }}$. Isomorphism classes of lifts of $Y_{\alpha}\to X_{K_f }$ to an unramified $\mathbb{Z}/4\mathbb{Z}$-cover are in bijection with a choice of square root of $(\beta -\alpha )$ for all $\beta \neq \alpha $, modulo the equivalence relation
\[
\{ (\beta -\alpha )^{1/2}\}\sim \{ - (\beta -\alpha )^{1/2} \}.
\]
\end{lemma}
\begin{proof}
Isomorphism classes of $\mathbb{Z}/2\mathbb{Z}$ covers of $Y_{\alpha ,K^{\sep }}$ unramified outside infinity are in bijection with equivalence classes of subsets of the set of roots of $\prod _{\beta \neq \alpha }(u_{\alpha }^2 -(\beta -\alpha )/c_{\alpha })$, under the equivalence relation $S_1 \sim S_2 $ if and only if $\Roots (g_{\alpha })=S_1 \sqcup S_2 $. This bijection respects the action of $\Aut (Y_{\alpha }/X_{K_f })$ on both sets. The $\mathbb{Z}/2\mathbb{Z}$-covers of $Y_{\alpha ,K^{\sep }}$ which come from a $\mathbb{Z}/4\mathbb{Z}$-cover of $X_{K^{\sep }}$ are exactly the $\Aut (Y_{\alpha }/X_{K_f })$-stable covers which do not descend to a $\mathbb{Z}/2\mathbb{Z}$-cover of $X_{K_f }$. The lemma follows from the fact that the generator of $\Aut (Y_{\alpha }/X_{K_f })$ acts on the set of roots of $\prod _{\beta \neq \alpha }(u_{\alpha }^2 -(\beta -\alpha )/c_{\alpha })$ by swapping $(\frac{\beta -\alpha }{c_{\alpha }})^{1/2}$ and $-(\frac{\beta -\alpha }{c_{\alpha }})^{1/2}$ for all $\beta \neq \alpha $.

\end{proof}
This construction, applied to all roots $\alpha $ of $f$, gives the function field of the maximal unramified $\mathbb{Z}/4\mathbb{Z}$-cover of $X_{K^{\sep }}$. Hence we deduce that the function field of the cover of $X_{K^{\sep }}$ corresponding to $J[4]$ is equal to 
\[
K^{\sep }(X)(u_{\alpha },\sqrt{g_{\alpha ,0}}:\alpha \in \Roots (f)).
\]

A modification of this construction also gives the cover of $X_K$ corresponding to $\Gamma $. Given distinct roots $\alpha $ and $\beta $ of $f$, we may consider the degree two extension of $K^{\sep }(X)(u_{\alpha },u_{\beta })$ generated by $z_{\alpha \beta }$ such that 
\[
z_{\alpha \beta }^2 :=c_{\alpha }c(u_{\alpha }-(\frac{\beta -\alpha }{c_{\alpha }})^{1/2}).
\]
Let $K^{\sep }(X)(u_{\alpha },u_{\beta },w_{\alpha \beta })$ be the field extension of $K^{\sep }(X)(u_{\alpha },u_{\beta })$ generated by $w_{\alpha \beta }$ satisfying 
\begin{equation}\label{eqn:vab}
w_{\alpha \beta }^2 :=2\left(1+ \frac{c_{\alpha }u_{\alpha }-c_{\beta }u_{\beta }}{\beta -\alpha}\right) .
\end{equation}
\begin{lemma}
We have an equality of function fields
\[
K^{\sep }(X)(u_{\alpha },u_{\beta },z_{\alpha \beta })=K^{\sep }(X)(u_{\alpha },u_{\beta },w_{\alpha \beta }).
\]
\end{lemma}
\begin{proof}
We have \begin{align*}
\Tr \left(\frac{u_{\alpha }-\sqrt{(\beta -\alpha )/c_{\alpha }}}{1-\sqrt{(\beta -\alpha )/c_{\alpha }}} \right) & =2\frac{c_{\alpha }u_{\alpha }+\alpha -\beta }{c_{\beta }} \\
\Nm \left(\frac{u_{\alpha }-\sqrt{(\beta -\alpha )/c_{\alpha }}}{1-\sqrt{(\beta -\alpha )/c_{\alpha }}} \right) & =u_\beta ^2 \\
\end{align*}
where the trace and norm are from $K(\alpha ,\beta ,\sqrt{(\beta -\alpha )/c_{\alpha }})(Y_{\alpha }\times _X Y_{\beta })$ to $K(\alpha ,\beta )(Y_{\alpha }\times _X Y_{\beta })$. 
Hence, by Lemma \ref{lemma:sqrt_nm}, the class of $\frac{u_{\alpha }-\sqrt{(\beta -\alpha )/c_{\alpha }}}{1-\sqrt{(\beta -\alpha )/c_{\alpha }}}$ in $K(\alpha ,\beta ,\sqrt{(\beta -\alpha )/c_{\alpha }})(Y_{\alpha }\times _X Y_{\beta })^\times \otimes \F_2$ is equal to that of 
\[
\frac{2(\alpha -\beta )}{c_{\beta }}\left(1+ \frac{c_{\alpha }u_{\alpha }-c_{\beta }u_{\beta }}{\alpha -\beta }\right).
\]
This is equal to $2\left(1+ \frac{c_{\alpha }u_{\alpha }-c_{\beta }u_{\beta }}{\alpha -\beta }\right)$ in $K^{\sep}(Y_{\alpha }\times _X Y_{\beta })^\times \otimes \F _2$. Finally, note that
\begin{align*}
& \left(1+ \frac{c_{\alpha }u_{\alpha }-c_{\beta }u_{\beta }}{\alpha -\beta }\right)\left(1+ \frac{c_{\alpha }u_{\alpha }-c_{\beta }u_{\beta }}{\beta -\alpha }\right) \\
= & -c_{\alpha }c_{\beta }\left( \frac{(u_{\beta }-1)(u_{\alpha }+1)}{(\alpha -\beta -c_{\beta} u_{\beta} -c_{\alpha }u_{\alpha })} \right) ^2 . 
\end{align*}
Hence 
\[
2\left(1+ \frac{c_{\alpha }u_{\alpha }(z)-c_{\beta }u_{\beta }(z)}{\alpha -\beta }\right) =\left(1+ \frac{c_{\alpha }u_{\alpha }(z)-c_{\beta }u_{\beta }(z)}{\beta -\alpha }\right)
\] in $K^{\sep}(Y_{\alpha }\times _X Y_{\beta })^\times \otimes \F _2$.
\end{proof}
Let $Y_{\alpha \beta }$ denote the degree $2$ cover of $X_{2,K_f ^{(2)}}$ defined by $w_{\alpha \beta }$. Let $Y_{\alpha ,2}$ denote the degree $2$ cover of $X_{2,K_f }$ defined by 
\[
h_{\alpha }^2 :=\Nm _{K_{f,2}(X_2 )|K_f (X_2 )} \left( 1+\frac{c_{\alpha }u_{\alpha }-c_{\beta }u_{\beta }}{\beta -\alpha } \right).
\]
Define
\[
X_4 :=\Res _{K_f |K}Y_{\alpha ,2}
\]
Then $X_4 $ is disconnected. Let $X_4 ^+ $ denote the subvariety of $X_4 $ defined by 
\[
\Nm _{K_f |K}(h_{\alpha })=\prod _{\{ \alpha ,\beta \} \in \Roots (f)^{(2)} }\left( 1+\frac{c_{\alpha }u_{\alpha }-c_{\beta }u_{\beta }}{\beta -\alpha } \right).
\]

Let $X_{\Gamma }$ be the cover of $X_{2,K^{\sep }}$ defined by 
\[
X_{\Gamma }:=\Res _{K_f ^{(2)}|K}Y_{\alpha \beta }\times _{\Res _{K_f ^{(2)}|K}X_2 }X_2 .
\]
Then $X_{\Gamma }$ is disconnected. Let $X_{\Gamma }^+$ denote the subvariety of $X_{\Gamma }$ defined by 
\[
\Nm _{K_f ^{(2)}(X_{\Gamma })|K_f (X_{\Gamma })}w_{\alpha \beta }=\prod _{\{\alpha ,\beta \} \in \Roots (f)^{(2)} }\left( 1+\frac{c_{\alpha }u_{\alpha }-c_{\beta }u_{\beta }}{\beta -\alpha } \right)
\]

Note that $Y_{\alpha }$ is isomorphic to the curve obtained by the pullback diagram
\begin{equation}
\begin{tikzcd}
Y_{\alpha } \arrow[d] \arrow[r] & \Res _{K_{f,2} |K_{f}}Y_{\alpha \beta ,K_{f,2}} \arrow[d, " 2( 1+\frac{c_{\alpha }u_{\alpha }-c_{\beta }u_{\beta }}{\beta -\alpha } ) " ] \\
\mathbb{P}^1 \arrow[r, " \Delta "] & \Res _{K_{f,2}|K}\mathbb{P}^1 .
\end{tikzcd}
\end{equation}

Then by construction we have a map $X_{\Gamma }' \to Y_{\alpha }$, where $X_{\Gamma }' /K_f$ is defined by the pullback diagram
\begin{equation}
\begin{tikzcd}
X_{\Gamma }' \arrow[d] \arrow[r] & \Res _{K_{f,2} |K_{f}}Y_{\alpha \beta ,K_{f,2}} \arrow[d] \\
X_2 \arrow[r, " \Delta "] & \Res _{K_{f,2}|K}X_{2,K_{f,2}} .
\end{tikzcd}
\end{equation}
The map $X_{\Gamma }' \to Y_{\alpha }$ is induced by the map 
\[
2^{2g}\Nm _{K_{f,2}(X)|K_f (X)}( 1+\frac{c_{\alpha }u_{\alpha }-c_{\beta }u_{\beta }}{\beta -\alpha } ): X_{\Gamma }' \to \mathbb{P}^1 ,
\] 
together with the map $X_{\Gamma }' \to \Res _{K_{f,2}|K_f }Y_{\alpha ,\beta ,K_{f,2}}$ above. Via the natural map $X_{\Gamma }\to \Res _{K_f |K}X_{\Gamma }'$, we obtain a map
\[
X_{\Gamma }\to X_4
\]
through which the covering $X_{\Gamma }\to X_2 $ factors.

We now show that $X_{\Gamma }\to X$ is the cover of $X$ corresponding to $\Gamma $. This can be reduced to a question on function fields, which can then be reduced to a group-theoretic question on Galois groups. In group-theoretic terms, we will reduce to the following lemma.
\begin{lemma}\label{lemma:FGamma}
Let $F$ be a free group on $2g$ generators $x_1 ,\ldots ,x_{2g}$. For each $i$, let $G_i$ denote the kernel of the homomorphism $F\to \mathbb{Z}/2\mathbb{Z}$ sending $x_i $ to $1$ and $x_j$ to $0$ for $j\neq i$. Let $H_i <G_i$ denote the kernel of the mod 2 abelianisation map 
\[
G_i \to G_i ^{\ab }\otimes \F _2 .
\]
Then 
\[
\Gamma \simeq \cap _i H_i .
\]
\end{lemma}
\begin{proof}
Clearly $H_i \subset \Ker (F\to \Gamma )$. Furthermore $\Ker (F\to F^{\ab }\otimes \Z /4\Z )\supset \cap _i H_i$, since the kernel of the surjection $F\to \mathbb{Z}/4\mathbb{Z}$ sending $x_i$ to $1$ and $x_j$ to 0 for $i\neq j$ contains $H_i$.

The kernel of $\Gamma \to F^{\ab }\otimes \Z /4\Z $ is an $\F _2 $-vector space with basis $[x_i ,x_j ]$ for $i\neq j$. On the other hand the image of $[x_i ,x_j ]$ in $H_i ^{\ab }\otimes \F _2 $ is nonzero. Hence $H_i \supset \Ker (F\to \Gamma )$.
\end{proof}

We now prove our result giving models for the (pointed) covers of $(X,b)$ corresponding to $\Gamma $ and $J[4]$.
\begin{lemma}\label{lemma:Gamma_descent}
\begin{enumerate}
\item A descent of the $\Gamma $-cover of $X_{K^{\sep }}$ to $K$ lifting the point $b$ is given by $X_{\Gamma }^+$.
\item A descent of the $J[4]$-cover of $X_{K^{\sep }}$ to $K$ lifting the point $b$ is given by $X_4 ^+ $.
\end{enumerate}
\end{lemma}
\begin{proof}
The proof is identical to that of Lemma \ref{lemma:X2}: we may base change to $K^{\sep }$ and apply Lemma \ref{lemma:WeilRes} to get an explicit description of $X_{\Gamma ,K^{\sep }}$ and $X_{4,K^{\sep }}$ as being given by the equations
\[
w_{\alpha \beta }^2 =2\left( 1+\frac{c_{\alpha }u_{\alpha }-c_{\beta }u_{\beta }}{\beta -\alpha } \right)
\]
for $\alpha ,\beta $ distinct roots of $f$, and
\[
h_{\alpha }^2 =\prod _{\beta \neq \alpha }\left( 1+\frac{c_{\alpha }u_{\alpha }-c_{\beta }u_{\beta }}{\beta -\alpha} \right)
\]
for $\alpha $ a root of $f$, respectively. We see from Lemma \ref{lemma:Z4Zcovers} that a generating set of unramified $\mathbb{Z}/4\mathbb{Z}$-covers of $X_{K^{\sep }}$ factor through . By Lemma \ref{lemma:FGamma}, the field generated by the function fields of the covers $Y_{\alpha \beta }$ is equal to that of the maximal $\Gamma $-cover of $X_{K^{\sep }}$.
\end{proof}

\subsection{Explicit formula for the $H^1 (K,\wedge ^2 J[2])$ class}
To complete the proof of Proposition \ref{prop:main2}, it remains to prove the formula for the lift of the nonabelian Kummer class $j_{\Gamma }(z)$ to $H^1 (K,\wedge ^2 J[2])$. By Lemma \ref{lemma:descent}, it is enough to construct a descent of the $\Gamma $-cover of $X_{K^{\sep }}$ to $K$ such that $b$ lifts to $K$-point, and describe the fibre.  The $\Gamma $-cover is given by $X_{\Gamma }^+$. To describe the fibre, first note that we have an identification of the preimage of $b' \in X_2 (K)$ with $\Ind ^{K_f ^{(2)}|K_f }\mu _2 $. Indeed by Lemma \ref{lemma:Gamma_descent} the fibre of $b'$ is isomorphic to the Weil restriction from $K_f ^{(2)}$ to $K$ of the fibre of $b'$ under $Y_{\alpha \beta }\to X_{K_f ^{(2)}}$. Hence the identification is a special case of the following lemma.

\begin{lemma}\label{lemma:shapiro_weil_restriction}
Let $L|K$ be a finite separable extension of fields, and $M$ a finite group scheme over $L$. 
\begin{enumerate}
\item 
We have an isomorphism of $\Gal (K^{\sep }|K)$-modules
\[
\Ind ^L _K (M(K^{\sep }))\simeq (\Res _{L|K}M)(K^{\sep }).
\]
\item
Suppose $P$ is an $M$-torsor over $L$. Give $\Res _{L|K}(P)$ the structure of a $\Res _{L|K}M$-torsor over $K$ via the map
\[
\Res _{L|K}(P) \times _K \Res _{L|K}(M)\to \Res _{L|K}(P)
\]
induced from the $M$-torsor structure on $L$. Then the class of $\Res _{L|K}(P)$ in $H^1 (K,\Ind ^L _K M(K^{\sep }))$ is the image of the class of $P$ in $H^1 (L,M(K^{\sep }))$ under the isomorphism from Shapiro's lemma.
\end{enumerate}
\end{lemma}
\begin{proof}
Part (1) follows from Lemma \ref{lemma:WeilRes} in its Galois equivariant form. For part (2), we can take the map $\Res _{L|K}(P)_L \to P$ adjoint to the identity morphism on $\Res _{L|K}(P)$. We obtain a $\Res _{L|K}(M)_L$-action on $P$ which we claim factors through the map $\Res _{L|K}(M)_L \to M$ adjoint to the identity, and recovers the original $M$-action on $P$. Indeed, base-changing to $K^{\sep }$ via an embedding $L\to K^{\sep }$, this can be seen from Lemma \ref{lemma:WeilRes}, since then (for any $S/L$ admitting a Weil restriction) the map $\Res _{L|K}(S)_{K^{\sep }}\to S_{K^{\sep }}$ is simply the map $\prod _{\sigma \in \Hom _K (L,K^{\sep })}S_{K^{\sep },\sigma }\to S_{K^{\sep }}$ projecting onto the factor corresponding to the chosen embedding $L\to K^{\sep }$. On the level of cocycles, the Lemma now follows from Lemma \ref{lemma:concrete_shapiro} below.
 \end{proof}

\begin{lemma}\label{lemma:concrete_shapiro}
Let $H<G$ be a finite index subgroup, and let $M$ be a finite $H$-module. Let $\pi :\Res ^G _H \Ind ^G _H M\to M$ be the map induced by adjunction.
The composite map
\[
H^1 (G,\Ind ^G _H M) \stackrel{\Res }{\longrightarrow }H^1 (H,\Res ^G _H \Ind ^G _H M) \stackrel{\pi _* }{\longrightarrow }H^1 (H,M) 
\]
is equal to the isomorphism from Shapiro's lemma.
\end{lemma}
\begin{proof}
See \cite{stix2010trading}, \cite{stix2013correction}, or \cite{NSW} for more general statements (for nonabelian cohomology and cohomology in arbitrary degree respectively).
\end{proof}

Above the point $b'$, the map $X_{\Gamma }^+ \to X_4$ is identified with the norm map
\[
\Ind ^{K_{f}^{(2)}}_{K }\F _2 \to \Ind ^{K_f }_K \F _2 .
\]
It follows that the the image of $z$ in $H^1 (K,\wedge ^2 J[2])$ is identified with the fibre of $z'\in X_2 (K)$ in $X_{\Gamma }(K)$, viewed as an $\Ind ^{K_f ^{(2)}}_K \F_2$-torsor whose class in $H^1 (K,J[2])$ vanishes. Hence Proposition \ref{prop:main2}
 follows from Lemma \ref{lemma:shapiro_weil_restriction} together with the description of the covers $X_{\Gamma }^+ \to X_2 $ and $X_4 ^+ \to X_2 $ above.

\section{$2$-descent for Selmer schemes: local aspects away from $2$}

Let $K$ be a number field, and $S$ a finite set of primes containing all primes above a rational prime $p$. Let $T$ be a free finite rank $\mathbb{Z}_p$-module with a continuous action of $G_{K,S}:=\Gal (K_S |K)$. Let $V:=T\otimes _{\Z _p }\Q _p $. We have short exact sequences of $G_{K ,S}$-modules
\begin{align}
0\to T\stackrel{\cdot p}{\longrightarrow }T\to T\otimes \F _p \to 0 \\
0\to T\to V\to V/T \to 0.
\end{align}
For any $G<G_{K ,S}$ and $i\geq 0$, we obtain an exact sequence
\begin{equation}\label{eqn:Hi}
0\to H^i (G,T)\otimes \F _p \to H^i (G,T\otimes \F _p )\to H^{i+1} (G,T)[p]\to 0
\end{equation}
and an isomorphism $H^i (G,T)\otimes _{\Z _p }\Q _p \simeq H^i (G,V)$. We want to bound the rank of $H^1 _f (G_{K ,S},V)$ by controlling the rank of $H^1 _f (G_{K ,S},T)\otimes \F _p $. Our goal is to find subspaces of $H^1 (G_{K ,S},T\otimes \F _p )$ which approximate $H^1 _f (G_{K ,S},T)\otimes \F _p $ and therefore whose $\F _p$-ranks give upper bounds on the $\Q _p $-rank of $H^1 _f (G_{K ,S},V)$. 

\begin{definition}
For $v$ a prime of $K$ not lying above $p$, and a $\Gal (K_v )$ module $M$, we define
\[
H^1 _{\ur}(K_v ,M):=H^1  (k_v ,M^{I_v })\subset H^1 (K_v ,M).
\]
For $T$ as above, we define
\[
H^1 _f (K_v ,T ):=i^{-1}H^1 _{\ur } (K_v ,V).
\]
\end{definition}
%
%
\begin{definition}
For $v$ a prime of $K$ not lying above $p$, we define $H^1 _{f,T} (K_v ,T\otimes \F _p )$ to be the image of $H^1 _f (K_v ,T)$ in $H^1 (K_v ,T/pT)$. We will sometimes abbreviate this to $H^1 _f (K_v ,T\otimes \F _p )$ (when the lattice $T$ is unambiguous) but it should be emphasised that this subspace (unlike the subspace $H^1 _{\ur}(K_v ,T\otimes \F_p )$) depends on $T$ and is not an invariant of the Galois module $T\otimes \F_p $ in general.
\end{definition}

\subsection{Properties of $H^1 _f (K_v ,T)$}
\begin{lemma}\label{lemma:unipotent_implies_same}
Let $T$ be as above. Let $v$ be a place of $K$, $k$ the residue field of the ring of integers of $K_v$, and $\sigma $ a topological generator for the tame inertia subgroup of the Galois group of $K_v$.
\begin{enumerate}
\item 
\begin{equation}
H^1 _f (K_v ,T\otimes \F_p )\subset \Ker (H^1 (K_v ,T\otimes \F _p )\to H^2 (K_v ,T)[p])
\end{equation}
\item Suppose $v$ is prime to $p$, $I_v $ acts unipotently on $T$, and $T/(\sigma -1)T$ is torsion-free. Then
\[
H^1 _f (K_v ,T\otimes \F_p ) \subset H^1 _{\ur}(K_v ,T\otimes \F _p ).
\]
\end{enumerate}
\end{lemma}
\begin{proof}
Part 1 is immediate from \eqref{eqn:Hi}. For part 2, by definition $H^1 _f (K_v ,T)$ consists of the classes whose image in $H^1 (I_v ,T)^{\Gal(\overline{k}_v |k_v )}$ is torsion. Hence if the latter is torsion-free, then $H^1 _f (K_v ,T)=H^1 _{\ur }(K_v ,T)$. Since the image of $H^1 _{\ur }(K_v ,T)$ in $H^1 (K_v ,T\otimes \F _p )$ is contained in $H^1 _{\ur}(K_v ,T\otimes \F _p )$, it is hence enough to show that $H^1 (I_v ,T)^{\Gal(\overline{k}_v |k_v )}$ is torsion-free when $I_v $ acts unipotently. Since $I_v$ acts unipotently, we have
\[
H^1 (I_v ,T)\simeq H^1 (I_v ^{(p)},T),
\]
where $I_v ^{(p)}$ is the maximal pro-p quotient of $I_v$. Hence we have an isomorphism of $\Z _p $-modules
\[
H^1 (I_v ,M)\simeq M/(\sigma -1)M
\]
where $\sigma $ is a topological generator of $I_v ^{(p)}$. Since $M/(\sigma -1)M$ is a torsion-free $\mathbb{Z}_p $-module, so are the invariants (by any group action).
\end{proof}

We will write $H^1 _f (K_v ,J[2])$ and $H^1 _f (K_v ,\wedge ^2 J[2])$ to mean the image of $H^1 _f (K_v ,T_2 J)$ and $H^1 _f (K_v ,\wedge ^2 T_2 J)$ respectively.

Recall that, given an oriented graph $\Gamma$ with vertex set $V$ and edge set $E$, we have maps
\begin{align*}
\delta :\Z [E]\to \Z [V] \\
e\mapsto s(e)-t(e).
\end{align*}
and 
\begin{align*}
\epsilon :\Z [V]\to \Z [E] \\
v \mapsto \sum _{s(e)=v}e-\sum _{t(e)=v}v
\end{align*}
We define 
\begin{align*}
H_1 (\Gamma ,\Z ):=\Ker (\delta ) \\
H^1 (\Gamma ,\Z ):=\Coker (\epsilon )
\end{align*}
and define the monodromy map
\[
N_{\Gamma }:H_1 (\Gamma ,\Z )\to H^1 (\Gamma ,\Z )
\]
to be the composite of the inclusion of $H_1 (\Gamma ,\Z )$ into $\Z [E]$ with the quotient map $\Z [E]\to H^1 (\Gamma ,\Z )$. Although the definitions above are phrased in terms of an oriented graph, for two different choices of an orientation on a graph $\Gamma $ (i.e. choices of an ordering of the endpoints of each edge) the homology and cohomology groups above are isomorphic in a manner compatible with the monodromy map. For this reason we will view $H_1 (\Gamma ,\Z ),H^1 (\Gamma ,\Z )$ and $N_{\Gamma }$ as having input a graph $\Gamma $.

Now let $K_v$ and $\sigma $ be as above. Let $\mathcal{O}_v$ denote the ring of integers of $K_v$ and $k _v$ the residue field. Let $X/K_v$ be a smooth projective geometrically irreducible curve with strict semistable reduction and $X/\mathcal{O}_v$ be a minimal regular model. Let $\Gamma $ denote the dual graph of $X_{k_v}$. Let $J$ be the Jacobian of $X$ and let $p$ be a prime which is a unit in $\mathcal{O}_v$. Grothendieck's Picard--Lefschetz formula describes the action of $\sigma $ on $T_p J$ in terms of $N_{\Gamma }$.
\begin{theorem}\label{thm:monodromy}[Picard--Lefschetz formula, \cite{SGA7}]

$T_p J$ admits a Galois stable filtration
\[
T_p J\supset W_1 T_p J \supset W_2 T_p J. 
\]
The map
\[
(\sigma -1):T_p J\to T_p J
\]
has kernel equal to $W_1 T_p J$ and image contained in $W_2 T_p J$. We have Galois equivariant isomorphisms
\begin{align*}
T_p J/W_1 T_p J\simeq H_1 (\Gamma ,\Z )\otimes _{\Z }\Z _p  \\
W_2 T_p J\simeq H^1 (\Gamma ,\Z )\otimes _{\Z }\Z _p (1)
\end{align*}
giving a commutative diagram
\[
\begin{tikzcd}
T_p J/W_1 T_p J \arrow[d, "\simeq "] \arrow[r, "\sigma -1"] \arrow[d] & W_2 T_p J \arrow[d, "\simeq"] \\
H_1 (\Gamma ,\Z )\otimes _{\Z }\Z _p \arrow[r, "N_\Gamma \otimes \alpha "]                                     & H^1 (\Gamma ,\Z )\otimes _{\Z }\Z _p (1)                    
\end{tikzcd}
\]
where $\alpha \in \Z _p (1)$ is an element corresponding to $\sigma $.
\end{theorem}
\begin{proof}
For the convenience of the reader we recall how to extract this statement from \cite[Expos\'e VIII Theorem 12.5]{SGA7}. First let $\mathcal{A}$ denote the Neron model of $J$ over $\mathcal{O}_v$. Let $\mathcal{A}^o _{k_v} \subset \mathcal{A}_{k_v} $ denote the connected component of the identity of the special fibre, and let $T\subset \mathcal{A}^o _{k_v}$ denote the toric part. Then the filtration is given by 
\begin{align*}
W_1 T_p J=T_p \mathcal{A}^o _{k_v} \\
W_2 T_p J=T_p T.
\end{align*}

In loc. cit., it is proved that two monodromy pairings coincide. The first is the bilinear pairing \cite[12.4.5]{SGA7}
\[
\Delta : H_1 (\Gamma ,\Z )\times H_1 (\Gamma , \Z )\to \Z 
\]
obtained by restriction from the bilinear pairing
\[
\Delta :\Z [E]\times \Z [E] \to \Z
\]
which sends $(e_1 ,e_2 )$ to $1$ if $e_1 =e_2 $ and $0$ otherwise.

On the other hand the Weil pairing defines an isomorphism $\Hom (W_2 T_p J ,\Z _p (1))\simeq T_p J/W_1 T_p J$, and hence $\sigma -1$ defines a map
\[
T_p J/W_1 T_p J\to \Hom (T_p J/W_1 T_p J ,\Z _p (1)),
\]
or dually a pairing 
\[
T_p J/W_1 T_p J \times T_p J/W_1 T_p J \to \Z _p (1).
\]
The statement of the Picard--Lefschetz theorem in loc. cit. implies that this pairing is identified with
\[
\Delta \otimes \alpha :H_1 (\Gamma ,\Z )\otimes \Z _p \times H_1 (\Gamma ,\Z )\otimes \Z _p \to \Z _p (1)
\]
via an isomorphism $H_1 (\Gamma ,\Z )\otimes \Z _p \simeq T_p J/W_1 T_p J$. Dualising gives the statement of the theorem above.
\end{proof}

\begin{corollary}
\begin{enumerate}
\item Suppose $X/K_v $ is a smooth projective geometrically irreducible curve with regular model whose special fibre is a curve of genus $g-1$ with two points identified. Let $T=\wedge ^2 T_p J$. Then $I_v$ acts unipotently on $T$, and $T/(\sigma -1)T$ is torsion free.
\item If $v$ lies above an odd prime, $f\in \mathcal{O}_v [x]$ is a polynomial of degree $\geq 5$ which is separable over $K_v$ with leading coefficient a $v$-adic unit, and $\val _v (\Delta (f))=1$, then $H^1 _f (K_v ,\wedge ^2 J[2])=H^1 _{\ur }(K_v ,\wedge ^2 J[2])$.
\end{enumerate}
\end{corollary}
\begin{proof}
For part (1), let $\Gamma $ denote the dual graph of a minimal regular model of $X$ over $\Z _{\ell }$. Then the monodromy map
\[
N:H_1 (\Gamma ,\Z  )\to H^1 (\Gamma ,\Z  )
\]
is an isomorphism. This implies $(\sigma -1)T_p J$ is torsion-free, and hence the same is true of $(\sigma -1)\wedge ^2 T_p J$. For part (2), note that the condition on the discriminant implies that $y^2 =f(x)$ defines a minimal regular model as in part (1). Hence the conclusion follows from part (1), Theorem \ref{thm:monodromy} and Lemma \ref{lemma:unipotent_implies_same}.
\end{proof}

\section{$2$-descent for Selmer schemes: local aspects at $2$}\label{sec:2}

To prove Theorem \ref{thm:BMSST_rational}, it will be necessary to describe $H^1 _f (\Q _2 ,\wedge ^2 J[2])$ as a subspace of $\Q _{2,f}^{(2),\times }\otimes \F _2 $. Recall that by definition $H^1 _f (\Q _2 ,\wedge ^2 J[2])$ is the image of $H^1 _f (\Q _2 ,\wedge ^2 T_2 J)$ in $H^1 (\Q _2 ,\wedge ^2 J[2])$. We aim to achieve this using the nonabelian $(x-T)$ map. The idea behind our approach is as follows. We can calculate the dimension of the image of $H^1 _f (\Q _2 ,\wedge ^2 T_2 J)$ in $H^1 (\Q _2 ,\wedge ^2 J[2])$, hence to describe the image it is enough to construct enough classes in $H^1 (\Q _2 ,\wedge ^2 J[2])$ which have crystalline lifts. We do this by using $\Q _2 $ points of $X$ to construct crystalline classes in a larger Galois cohomology set, and try to give conditions for these to `induce' crystalline classes with values in $\wedge ^2 J[2]$.

More precisely, let $b$ and $z$ be $\Z _2 $ points of $X$ not reducing to the point at infinity, and let $\Pi _2 $ denote the maximal $2$-unipotent quotient of $\pi _1 ^{\et ,2}(X_{\overline{\Q }_2 }-\{ \infty \} ,b)$. Then $j_{\Pi _2 ,b}(z)$ is an element of $H^1 _f (\Q _2 ,\Pi _2 )$, i.e. an element of $H^1 (\Q _2 ,\Pi _2 )$ whose image in $H^1 (\Q _2 ,U_2 )$ lies in the subset $H^1 _f (\Q _2 ,U_2 )$ in the sense of Kim \cite{kim:chabauty}. Here $U_2$ denote the maximal $2$-nilpotent quotient of the $\Q _p $-unipotent fundamental group of $X_{\overline{\Q }_2 }-\infty $. Furthermore, if $z-b$ lies in $4\cdot J(\Q _2 )$, then $j_{\Gamma ,b}(z)$ lifts to an element of $H^1 (\Q _2 ,\wedge ^2 J[2])$, which we shall also denote $j_{\Gamma ,b}(z)$ (although if $J[4](\Q _2 )\neq 0$ then this lift is not unique).

Now suppose we have another class $\gamma \in H^1 _f (G_{\Q _2 },\Pi _2 )$ such that the image of $\gamma $ in $H^1 (\Q _2 ,T_2 J)$ is equal to that of $j_{\Pi _2 ,b}(z)$. Then we may twist $\Pi _2 $ by a cocycle defining $\gamma$ to get a new group $\Pi _2 ^{(\gamma )}$ which is still a central extension of $T_2 J$ by $\wedge ^2 T_2 J$, and under the bijection between $H^1 (\Q _2 ,\Pi _2 )$ and $H^1 (\Q _2 ,\Pi _2 ^{(\gamma )})$ (see e.g. \cite[I.5.34]{serre-gc}), $j_{\Gamma ,b}(z)$ will be sent to a new element which we denote $\gamma ^{-1}j_{\Gamma ,b}(z)$. Then $\gamma ^{-1}j_{\Gamma ,b}(z)$ defines an element of $H^1 _f (\Q _2 ,\wedge ^2 T_2 J)$. Now suppose furthermore that the image of $\gamma $ in $H^1 (\Q _2 ,\Gamma )$ is trivial. Then the image of $\gamma ^{-1}j_{\Gamma ,b}(z)$ in $H^1 (\Q _2 ,\wedge ^2 J[2])$ is equal to that of $j_{\Gamma ,b}(z)$. Hence we have verified the existence of a crystalline lift of $j_{\Gamma ,b}(z)$.
\subsection{Integral mixed extensions}
The issue with the idea outlined above is that the only way we have of constructing elements of $H^1 _{f} (\Q _2 ,\Pi _2 )$ is from $\Z _2 $-points of the curve $X$. This means that the condition of finding a $\gamma $ with the same image in $H^1 (\Q _2 ,T_2 J)$ as $[z]-[b]$ is too restrictive. For this reason, we replace $\Pi  _2 $ with a larger group for which it is easier to construct elements. 

Let $\Lambda =\F _p ,\Z _p $ or $\Q _p $.
Let $T_0 ,\ldots ,T_n $ be finite rank free $\Lambda $ module with a continuous action of a profinite group $G$. We define a mixed extension with graded pieces $T_0 ,\ldots ,T_n $ to a be a tuple $((M_\bullet ), (\psi _\bullet ))$ consisting of a finite rank free $\Lambda $-module with a continuous action of $G$ and a $G$-stable filtration
\[
M_0 \supset M_1 \supset \ldots \supset M_n = 0
\]
together with $G$-equivariant isomorphisms of $\Lambda $-modules $\psi _i :M_{i-1}/M_i \stackrel{\simeq }{\longrightarrow }T_i$ for all $1\leq i\leq n$. Given two mixed extensions $(M_\bullet ,\psi _\bullet ),(M' _\bullet ,\psi ' _\bullet )$ with graded pieces $T_0 ,\ldots ,T_n $, we define a unipotent isomorphism between $M$ and $M'$ to be an isomorphism of filtered $\Lambda $-modules $f:M_0 \simeq M_0 '$ which is compatible with the $\psi _i $ and $\psi _i '$. By compatible we mean that, if $\gr _i f$ denotes the induced isomorphism $M_{i-1}/M_i \simeq M_{i-1}' /M_i '$, then we have $\psi _i =\psi _i ' \circ \gr _i f$ for all $i$.

Given a mixed extension $(M_\bullet ,\psi _\bullet )$ as above, we define the group $U(M_\bullet )$ to be the group of unipotent isomorphisms from $M_\bullet $ to itself. When $M_\bullet $ is the trivial mixed extension $\oplus T_i$, we denote this as $U(T_\bullet )$. The groups $U(M_\bullet )$ naturally have $G$-actions, and for two mixed extensions $M_\bullet $ and $M_\bullet '$, the set of unipotent isomorphisms from $M_\bullet $ to $M_\bullet '$ naturally has the structure of a Galois-equivariant $(U(M_\bullet ),U(M_\bullet '))$-bitorsor (or $(U(M_{\bullet} ),U(M_{\bullet }'))$-principal space in the sense of \cite[I.5.35]{serre-gc}). The following lemma may be proved as in \cite[Lemma 4.7]{QC1}.
\begin{lemma}
For any mixed extension $M_\bullet $,we have a bijection between the set of isomorphism classes of mixed extensions with graded pieces $T_\bullet $ and the nonabelian cohomology set $H^1 (G,U(M_\bullet ))$, defined by sending a mixed extension $M_\bullet '$ to the $G$-equivariant $U(M_\bullet )$-torsor of unipotent isomorphisms between $M_\bullet $ and $M_\bullet '$.
\end{lemma}
Now let $\Pi _n $ be the maximal pro-$p$ $n$-unipotent quotient of $\pi _1 ^{\et }(X_{\overline{K}},b)$. Let $T_i$ denote the $i$th graded piece of the Iwasawa algebra $\Z _p [\! [\pi _1 ^{\et }(X_{\overline{K}}-\{ \infty \} ,b)]\! ]$ with respect to the $I$-adic filtration, where $I$ is the augmentation ideal. Define 
\[
E_n (b):=\Z _p [\! [\Pi _n ]\! ]/I^{n+1}.
\]
Then $E_n (b)$, via its $I$-adic filtration, acquires the structure of a mixed extension with associated graded $T_\bullet $. Via the unipotent action of $\Pi _n$ on $E_n (b)$, we obtain a $G$-equivariant homomorphism
\[
\Pi _n \to U(E_n (b)).
\]
Hence we have a map
\[
H^1 (G,\Pi _n )\to H^1 (G,U(E_n (b))).
\]
The image of the class of $\pi _1 ^{\et }(X_{\overline{K}}-\{ \infty \};b,z)$ is the class of the mixed extension
\[
E_n (b,z):=E_n (b)\times _{\pi _1 ^{\et }(X_{\overline{K}}-\{ \infty \},b)}\pi _1 ^{\et }(X_{\overline{K}}-\{ \infty \};b,z).
\]
The tensor product $E_n (b)\otimes _{\Z _p }\F _p $ has the structure of a mixed extension with associated graded $T_\bullet \otimes _{\Z _p }\F _p $, and we have a $G$-equivariant surjection $U(E_n (b))\to U(E_n (b)\otimes \F _p )$. We define $\Gamma _n $ to be the image of $\Pi _n $ in $U(E_n (b)\otimes \F _p )$. For example, when $n=p=2$, $\Gamma _2 =\Gamma $. When $n=2$ and $p>2$, $\Gamma _2 $ will be an extension of $J[p]$ by $\wedge ^2 J[p]$.

Now suppose $n=2$. Then $U(T_\bullet )$ is a central extension
\[
1\to T_0 ^* \otimes T_2 \to U(T_\bullet )\to T_0 ^* \otimes T_1 \oplus T_1 ^* \otimes T_2 \to 1.
\]
In terms of mixed extensions, a lift of an element in $H^1 (G,T_0 ^* \otimes T_1 )\times H^1 (G,T_1 ^* \otimes T_2 )$ to $H^1 (G,U(T_{\bullet })$ corresponds to concatenating an extension $[E_1 ]$ of $T_0 $ by $T_1$ and an extension $[E_2 ]$ of $T_1 $ by $T_2 $ to form a mixed extension with associated graded $T_\bullet $. Following Nekov\'{a}\v{r} \cite{nekovar} (following SGA \cite{SGA7}), we shall refer to such a concatenation as a mixed extension of $E_1 $ and $E_2$. The obstruction to forming such a mixed extension is exactly the image of $([E_1 ],[E_2 ])$ in $H^2 (G,T_0 ^* \otimes T_2 )$ under the boundary map. We recall that this boundary map has the following simple description.
\begin{lemma}
The boundary map
\[
H^1 (G,T_0 ^* \otimes T_1 \oplus T_1 ^* \otimes T_2 )\to H^2 (G,T_0 ^* \otimes T_2 )
\]
is given by 
\[
(c_1 ,c_2 )\mapsto c_1 \cup c_2 ,
\]
In particular, a pair $(c_1 ,c_2 )\in H^1 (G,T_0 ^* \otimes T_1 \oplus T_1 ^* \otimes T_2 )$ lies in the image of $H^1 (G,U(T_0 ,T_1 ,T_2 ))$ if and only if $c_1 \cup c_2 $ is zero in $H^2 (G,T_0 ^* \otimes T_2 )$ with respect to the pairing defined by the natural map
\[
T_0 ^* \otimes T_1 \otimes T_1 ^* \otimes T_2 \to T_0 ^* \otimes T_2 .
\]
\end{lemma}
\begin{proof}
By definition, in terms of inhomogeneous cochains, the boundary map can be given by first choosing a (set-theoretic) section of 
\[
U(T_0 ,T_1 ,T_2 )\to T_0 ^* \otimes T_1 \oplus T_1 ^* \otimes T_2 
\]
which we take to be, in matrix notation
\[
s:(\alpha ,\beta )\mapsto \left( \begin{array}{ccc}1 & 0 & 0 \\ \alpha & 1 & 0 \\ 0 & \beta & 1 \end{array} \right).
\]
Then the boundary map sends $(c_1 ,c_2 )$ to the $2$-cocycle
\[
(g,h)\mapsto \beta (g)\cdot (g\cdot \alpha (h)),
\]
which is exactly the cup product.
\end{proof}
When working with mixed extensions over $\Q _p $, one can rescale the isomorphisms identifying the associated graded with $\oplus T_i$ to obtain a new mixed extension structure on the same representation. One new subtlety in the case of mixed extensions with $\Z _p$ rather than $\Q _p$ coefficients is that, given a mixed extension $M$, one can also construct new mixed extensions by choosing different lattices inside $M\otimes \Q _p $. This is used later to kill torsion. A simple case of this idea is presented in the following lemma.
\begin{lemma}\label{lemma:trivial}
If $M$ is an integral crystalline mixed extension with graded pieces $T_0 ,T_1 ,T_2 $, such that $[M/T_2 ]\in p\cdot H^1 _f (\Q _p ,T_0 ^* \otimes T_1 )$, then there is a Galois stable sublattice $M'$ of $M$ with $M/M' \simeq T_1 \otimes \F _p $, and such that
\[
p[M' /T_2 ]=[M/T_2 ]
\]
in $H^1 _f (\Q _p ,T_0 ^* \otimes T_1 )$ and 
\[
[\Ker (M'\to T_0 )]=p[\Ker (M\to T_0 )]
\]
in $H^1 _f (\Q _p ,T_1 ^* \otimes T_2 )$.
\end{lemma}
\begin{proof}
One may choose a $\Z _p $-module splitting $\alpha :T_0 \oplus T_1 \stackrel{\simeq }{\longrightarrow }M/T_2 \simeq T_0 \oplus T_1$ with respect to which the action of $G_{\Q _p }$ is given by
\[
\left( \begin{array}{cc} \rho _0 & 0 \\ p\cdot c(g) & \rho _1 \end{array} \right) .
\] 
Then $\alpha (T_0 \oplus pT_1 )$ is a Galois stable sublattice. Taking the preimage of $\alpha (T_0 \oplus pT_1 )$ in $M$ gives the lemma.
\end{proof}

We now turn to the problem of describing the structure of $E_2 (b)$. By definition this is a direct sum of $\Z _2 $ and $IE_2 (b)$. Hence the problem is to describe the structure of $I\Z _2 [\! [\Pi _2 ]\! ]/I^3 $ as an extension of $T_2 J$ by $T_2 J^{\otimes 2}$. This is closely related to the class of the Ceresa cycle of $(X,b)$ in Galois cohomology, by work of Hain and Matsumoto \cite{HM}. The relation \textit{integrally} (rather than over $\Q _2 $) has some additional subtleties which were recently studied by Bisogno, Li, Litt and Srinivasan \cite{BLLS}. They show that, if $\infty $ denotes a rational Weierstrass point, then the class of $I\Z _2 [\! [\Pi _2 (\infty ) ]\! ]/I^3 $ is $2$-torsion.

Let $\theta :H^1 _f (\Q _2 ,T_2 J)\times H^1 _f (\Q _2 ,T_2 J)\to H^1 _f (\Q _2 ,\Hom (T_2 J,T_2 J^{\otimes 2}))$ 
be the map induced by the left and right multiplication maps
\begin{equation}\label{eqn}
T_2 J\times T_2 J\to \Hom (I/I^2 ,I^2 /I^3 ).
\end{equation}
We deduce that, for arbitrary $b,z\in (X-\infty )(\Z _2 )$, the class of $IE_2 (b,z)$ in $H^1 _f (\Q _2 ,\Hom (T_2 J,T_J^{\otimes 2}))$ is equal to 
\begin{equation}
[E_2 (\infty )]+\theta (z-\infty ,b-\infty ).
\end{equation}
In particular, this implies the following.
\begin{lemma}
For any $b,z\in (X-\infty )(\Z _2 )$, we have
\[
2[IE_2 (b,z)]\in 2\cdot \mathrm{Im}(\theta ).
\]
\end{lemma}

Suppose that $f$ is irreducible over $\Q _2$. Let $i _f :D_f \hookrightarrow  H^1 _f (\Q _2 ,\Hom (T_2 J,I^2 /I^3 ))$ denote the image of $\theta $, and let $i : D\hookrightarrow H^1 (\Q _2 ,\Hom (T_2 J,I^2 /I^3 ))$ denote the image of $H^1  (\Q _2 ,T_2 J)^{\oplus 2}$. Let $E_f (T_2 J)$ be a crystalline extension
\[
0\to I^2 /I^3 \to E(T_2 J)\to T_2 J\otimes D_f \to 0
\]
such that the class of $E_f (T_2 J)$ in $H^1 _f (\Q _2 , \Hom (T_2 J\otimes D_f ,I^2 /I^3 ))$ is equal to $i _f$ under the isomorphism of $\Z _2 $-modules
\begin{align*}
& H^1 _f (\Q _2 ,\Hom (T_2 J\otimes D_f ,I^2 /I^3 )) \\
\simeq & \Hom _{\Z _2 } (D_f ,H^1 _f (\Q _2 ,\Hom (T_2 J,I^2 /I^3 ))).
\end{align*}
Similarly let $E(T_2 J)$ be an extension 
\[
0\to I^2 /I^3 \to E(T_2 J)\to T_2 J\otimes D \to 0
\]
such that the class of $E(T_2 J)$ in $H^1 (\Q _2 ,\Hom (T_2 J\otimes D,I^2 /I^3 ))$ is equal to $i$ under a similar isomorphism.

For any extension $E$ of $T_2 J$ by $I^2 /I^3 $ whose class in $\ext ^1 (T_2 J,I^2 /I^3 )$ lies in the image of $D$ we obtain a Galois equivariant map
\[
\iota _E  :E\to E(T_2 J).
\]
For any crystalline extension $E$ of $T_2 J$ by $I^2 /I^3 $ whose class in $\ext ^1 _f (T_2 J,I^2 /I^3 )$ lies in the image of $D_f$, we obtain a Galois equivariant map
\[
\iota _{E,f}:E\to E_f (T_2 J).
\]
These two maps are compatible via the natural injection
\[
\iota _f :E_f (T_2 J)\hookrightarrow E(T_2 J),
\]
 i.e. $\iota _E = \iota _f \circ \iota _{E,f}$. Given a mixed extension $M$ with graded pieces $\Z _2 ,T_2 J $ and $I^2 /I^3$, we obtain an extension
\[
0\to M_1 \to M\to \Z _2 \to 0
\]
defining an element of $H^1 (\Q _2 ,M_1 )$. If the class of $M_1 $ in $\ext ^1 (T_2 J,I^2 /I^3 )$ lies in the image of $D$ then we will denote by $[M]\in H^1 (\Q _2 ,E(T_2 J))$ the image of this extension class under the map
\[
\iota _{M_1 ,*}:H^1 (\Q _2 ,M_1 )\to H^1 (\Q _2 ,E(T_2 J)).
\]
If $M$ is moreover crystalline, we can play the same game instead using the map $\iota _{E,f}$ to obtain a class in $H^1 _f (\Q _2 ,E_f (T_2 J))$. We deduce the following characterisation of $H^1 _f (\Q _2 ,E_f (T_2 J))$.
\begin{lemma}\label{lemma:cris1}
An integral mixed extension $[M]$ whose class in $\ext ^1 (T_2 J,I^2 /I^3 )$ lies in the image of $D$ is crystalline if and only if its class in $H^1 (\Q _2 ,E (T_2 J))$ is in the image of $H^1 _f (\Q _2 ,E_f (T_2 J))$.
\end{lemma}

This implies the following criterion for the existence of crystalline mixed extensions lifting crystalline classes in $H^1 (\Q _2 ,T_2 J)$ and $D\subset H^1 (\Q _2 ,T_2 J ^* \otimes I^2 /I^3 )$.
\begin{lemma}\label{lemma:cris2}
If $c_1 \in H^1 _f (\Q _2 ,T_2 J )$ and $c_2 \in D_f$ have the property that there are crystalline mixed extensions $M_1 ,\ldots ,M_n$ such that 
\[
c_1 \otimes c_2 =\sum _i \pi _{1 *}[M_i ]\otimes \pi _{2 *}[M_i ],
\]
then $(c_1 ,c_2 )$ admit a lift to a crystalline mixed extension.
\end{lemma}
\begin{proof}
Since $c_1 \cup c_2 =\sum \pi _{1 *}(c_1 )\cup \pi _{2 *}(c_2 )=0$, there is an integral mixed extension $M$ of $c_1 $ and $c_2 $. Then $[M]-\sum [M_i ]$ is an element of $H^1 (\Q _2 ,E(T_2 J ))$ which lies in the image of $H^1 (\Q _2 ,I^2 /I^3  )$, by the exact sequence
\[
H^1 (\Q _2 ,I^2 /I^3  )\to H^1 (\Q _2 ,E(T_2 J))\to H^1 (\Q _2 ,T_2 J )\otimes D.
\]
In particular, we may twist $M$ by a class $\xi $ in $H^1 (\Q _2 ,I^2 /I^3 )$ to obtain a mixed extension $M^\xi $ with the property that $[M]=\sum [M_i ]$ in $H^1 (\Q ,E(T_2 J))$. In particular, since all the $M_i$ are in the image of $H^1 _f (\Q _2 ,E_f (T_2 J))$, so is $M^{\xi }$, and hence $M^{\xi }$ is crystalline by Lemma \ref{lemma:cris1}.
\end{proof}
%
%

%

We obtain the following criterion for liftability which is amenable to computations. Let $\widehat{J}$ denote the $2$-adic completion of $J(\Q _2 )$.

\begin{lemma}\label{lemma:cris_condition}
Suppose the homomorphism of $\Z _2 $-modules
\[
\Z _2 [(X-\infty )(\Z _2 )\times (X-\infty )(\Z _2 )]\to \widehat{J}^{\otimes 2} \times \widehat{J}^{\otimes 2}
\]
sending $(z_1 ,z_2 )$ to $((z_1 -z_2 )\otimes (z_1 -\infty ),(z_1 -z_2 )\otimes (z_2 -\infty ))$ has image equal to the preimage of $\Sym ^2 \widehat{J}\subset \widehat{J}^{\otimes 2}$ under the map
\[
\widehat{J}^{\otimes 2}\times \widehat{J}^{\otimes 2} \to \widehat{J}^{\otimes 2}
\]
sending $(v,w)$ to $v-w$. Then for ever $b,z\in (X-\infty )(\Z _2 )$ such that $z-b\in 4J(\Q _2 )$, the lift of $j_{\Gamma ,b}(z)$ to $H^1 (\Q _2 ,\wedge ^2 J[2])$ is in the image of $H^1 _f (\Q _2 ,\wedge ^2 T_2 J)$.
\end{lemma}
\begin{proof}
By \eqref{eqn} and Lemma \ref{lemma:cris2}, the conditions imply that there is a crystalline mixed extension $M$ who image in $H^1 _f (\Q _2 ,T_2 J)\otimes H^1 _f (\Q _2 ,T_2 J ^* \otimes T_2 J^{\otimes 2})$ is equal to twice that of $E_2 (b,z)$. Rescaling as in Lemma \ref{lemma:trivial} we may replace $E_2 (b,z)$ with a mixed extension $M'\subset E_2 (b,z)$ whose class in $H^1 _f (\Q _2 ,T_2 J^* \otimes I^2 /I^3 )$ lies in the image of $\theta $. Hence taking $[M']-2[M]$ gives an extension of $\Z _2 $ by $T_2 J^{\otimes 2}$ realising the extension class.
\end{proof}

\section{The Chabauty--Coleman--Kim method}

\subsection{Galois cohomology of number fields}
We first recall the conjectural description of the Galois cohomology of global Galois representations coming from geometry \cite{BK}.
\begin{conjecture}[Bloch--Kato]\label{conj:BK}
Let $S_0$ be a finite set of primes, and suppose $\mathcal{X}$ is a proper regular model over $\mathbb{Z}_{S_0 }$ of a smooth projective variety $X/\Q $. Then the \'etale regulator
\[
K_{2r-m-1} ^{(r)}(\mathcal{X})\otimes \Q _p \stackrel{\simeq }{\longrightarrow }H^1 _{f,S_0 } (G_\Q ,H^m (X_{\overline{\Q }},\Q _p (r)))
\]
is an isomorphism.
\end{conjecture}
Although we will mostly be interested in Galois cohomology for Galois representations of negative weight, the following form of Poitou--Tate duality will allow us to restrict to the above conjecture when the weight is non-negative, i.e. when $m\geq 2r$. Note that in this case the left-hand side is zero by definition, and hence the conjecture says that if $m\geq 2r$, then
\[
H^1 _{f,S_0 }(G_\Q ,H^m (X_{\overline{\Q }},\Q _p (r)))=0.
\]
\begin{proposition}
Let $W$ be a finite dimensional $\Q _p$ representation of $\Gal (K)$ unramified outside a finite set of primes.
\begin{enumerate}
\item \cite[Remark 1.2.4]{FPR94}

If $K=\Q $, and $H^0 (G_{\Q ,T},W)$, $H^0 (G_{\Q ,T},W^* (1))$ and $D_{\cris }(W)^{\phi =1}$ are all zero, then
\begin{align*}
& \dim H ^1 _f (\Q _p ,W)- \dim H ^1 _f (G_{\Q ,T},W) \\
= & \dim H^0 (\mathbb{R},W)- \dim H^1 _f (G_{\Q ,T},W^* (1)).
\end{align*}
\item \cite[1.2.2]{FPR94}

For any number field $K$ and $\Q _p $-Galois representation $W$,
\begin{align*}
& \dim \Ker (H^1 (G_{K,T},W)\to \oplus _{v\in T}H^1 (K_v ,W))  \\
= & \dim \Ker (H^2 (G_{K,T},W^* (1))\to \oplus _{v\in T}H^2 (K_v ,W^* (1))).
\end{align*}
\end{enumerate}
\end{proposition}

Hence we arrive at the following form of (part of) the Bloch--Kato conjectures.
\begin{conjecture}
For $X$ as above, $p$ a prime of good reduction and $m< 2r-1$,
\begin{align*}
& \dim H^1 _f (\Q ,H^m _{\et }(X_{\overline{\Q }},\Q _p (r)) \\
= & \dim H^1 _f (\Q _p ,H^m _{\et }(X_{\overline{\Q }},\Q _p (r)))- \dim H^0 (\R ,H^m (X_{\overline{\Q }},\Q _p (r))) \\
& -\dim H^0 (\Q ,H^m (X_{\overline{\Q }},\Q _p (r-1))^* ).
\end{align*}
\end{conjecture}
\subsection{The Chabauty--Coleman--Kim sets $X(\Q _p )_n $}

Let $X$ be a smooth projective geometrically irreducible curve over $\Q $, of genus $g>1$. Let $b\in X(\Q )$. Fix a prime $p$ of good reduction for $X$. The Chabauty--Coleman--Kim method, as developed in \cite{kim:chabauty}, \cite{kim2008component} and \cite{kim:siegel}, determines subsets $X(\Q _p )_n \subset X(\Q _p )$, containing the set of rational points $X(\Q )$, which are defined in terms of the $\Q _p $-unipotent fundamental group of $X$, which we now define.

For $n>0$, let $U_n (b)$ denote the maximal $n$-unipotent quotient of the $\Q _p $-unipotent fundamental group of $X$. This may be defined by the following universal property: it is a unipotent group together with a continuous homomorphism
\[
\rho _n :\pi _1 ^{\et }(X_{\overline{\Q }},b)\to U_n (b)(\Q _p ),
\]
such that for any $n$-unipotent group $U$ over $\Q _p $ with a continuous homomorphism 
\[
\pi _1 ^{\et }(X_{\overline{\Q }},b) \to U(\Q _p ),
\]
there is a unique homomorphism $U_n (b)\to U$ making 
\[
\begin{tikzcd}
\pi _1 ^{\et }(X_{\overline{\Q }},b) \arrow[r] \arrow[rd] & U_n (b) \arrow[d] \\
                       & U(\Q _p )                  
\end{tikzcd}
\]
commute. By universal properties, $U_n (b)$ inherits an action of $\Gal (\overline{\Q }|\Q )$ from $\pi _1 ^{\et }(X_{\overline{\Q }},b)$ making the homomorphism $\rho _n$ Galois-equivariant. Hence we obtain morphisms
\[
j_{n,K}:X(\Q )\to H^1 (\Q ,U_n (b)).
\]
Let $C_i U_n (b)$ denote the central series filtration, defined so that 
\[
U_n (b)/C_i U_n (b)\simeq U_i (b)
\]
for $i\leq n$. Finally, for $a\leq b\leq m\leq n$, let $s_{a,b,m,n}$ denote the homomorphism
\[
C_b U_n \to C_a U_m
\]
\begin{theorem}[Kim, \cite{kim:siegel}]
\begin{enumerate}
\item For $G=G_{\Q ,S}$ or $\Gal (\overline{\Q } _v |\Q _v )$ for a prime $v$ of $\Q $, $H^1 (G,U_n (b))$ is isomorphic to the $\Q _p$-points of a scheme of finite type over $\Q _p$, in such a way that all the maps $\loc _v$ and $s_{a,b,m,n*}$ are induced from morphisms of schemes.
\item $H^1 _f (G_{\Q _v },U_n (b))$ is (the $\Q _p $-points of) a $\Q _p$-subvariety of $H^1 (\Q _v ,U_n (b))$, and moreover this $\Q _p$-subvariety is isomorphic to an affine space.
\end{enumerate}
\end{theorem}
We define $\Sel (U_n (b))\subset H^1 (G_{\Q ,S},U_n (b))$ to be the fibre product, over $H^1 (G_{\Q ,S},U_1 (b))\times \prod _{v\in S}H^1 (\Q _v ,U_n (b))$, of $H^1 (G_{\Q ,S},U_n (b))$ with the following subschemes:
\begin{align*}
H^1 _f (\Q _p ,U_n (b))\subset H^1 (\Q _p ,U_n (b)), \\
j_{n,\Q _v }(X(\Q _v ))\subset H^1 (\Q _v ,U_n (b)), & v\in S-\{p\} \\
\mathrm{Im}(\Jac (X)(\Q )\otimes \Q _p )\subset H^1 (G_{\Q ,S},U_1 (b)). &  \\
\end{align*}
We define $X(\Q _p )_n$ to be the set of points of $X(\Q _p )$ whose image in $U_n ^{\dR}/F^0 $ lies in the image of $\Sel (U_n )$. By construction $X(\Q _p )_n$ contains $X(\Q )$.
\begin{theorem}[\cite{kim:chabauty}]
Suppose 
\[
\rk J(\Q ) +\sum _{i=2}^n \dim H^1 _f (\Q  ,\gr _i U_n ) <\sum _{i=1}^n \dim H^1 _f (\Q _p ,\gr _i U_n ).
\]
Then $X(\Q _p )_n$ is finite.
\end{theorem}
The right hand side equals $\dim U_n -\dim F^0 U_n $. We know that $U_n$ is isomorphic to the maximal $n$-unipotent quotient of a free unipotent group on $2g$ generators modulo one relation in depth 2, and $F^0 U_n$ is isomorphic to the maximal $n$-unipotent quotient of a free unipotent group on $g$ generators.

It follows that, when $n=2$, the dimension of the right hand side is
\[
\frac{(3g-2)(g+1)}{2}.
\]
We deduce the following.
\begin{lemma}
The set $X(\Q _p )_2 $ is finite whenever
\[
\dim H^1 _f (\Q ,\wedge ^2 V)<\frac{(3g-2)(g+1)}{2}-\rk J(\Q ).
\]
\end{lemma}
\subsection{Dimension bounds via descent}
We have the following result (see \cite[Proposition 12.6]{PS97}).
\begin{proposition}\label{prop:PS}
Let $L$ be a number field. For any finite set $S$ of primes containing all archimedean places, we have a short exact sequence
\[
0\to \mathcal{O}_{L,S}^\times /\mathcal{O}_{L,S}^{\times p} \to (L^\times /L^{\times p})_S \to \Cl (\mathcal{O}_{L,S})[p]\to 0.
\]
\end{proposition}
In our case of interest this gives the following dimension formula.
\begin{lemma}\label{lemma:PS2}
Let $J$ be the Jacobian of a hyperelliptic curve over a number field $K$ with an odd degree model of the form $y^2 =f(x)$, and Jacobian $J$. 
We have
\begin{align*}
\dim H^1 (\mathcal{O}_{K,S},\wedge ^2 J[2]) & =\rk \mathcal{O}_{K_f ^{(2)},S}^\times +\rk \mathcal{O}_{K,S}^\times -\rk \mathcal{O}_{K_f ,S}^\times \\
& + \dim \Cl (\mathcal{O}_{K_f ^{(2)},S})[2]-\dim \Cl (\mathcal{O}_{K_f ,S})[2]+\dim \Cl (\mathcal{O}_{K,S})[2].
\end{align*}
\end{lemma}
\begin{proof}
From the direct sum decompositions $\Ind ^{K_f }_K \mu _2 \simeq J[2]\oplus \mu _2 $ and $\Ind ^{K_f ^{(2)}}_K \mu _2 \simeq \wedge ^2 J[2]\oplus J[2]$ we obtain 
\[
\dim H^1 (\mathcal{O}_{K,S},\wedge ^2 J[2])=\dim H^1 (\mathcal{O}_{K_f ^{(2)}},\mu _2 )+\dim H^1 (\mathcal{O}_{K},\mu _2 )-\dim H^1 (\mathcal{O}_{K_f },\mu _2 ).
\]
The lemma then follows from Proposition \ref{prop:PS}.
%
\end{proof}

\subsection{Generalised heights}
We briefly describe how to explicitly compute $X(\Q _2 )_2$ using the above. More details may be found in \cite{BKC}. Suppose for simplicitly that $X$ is a smooth hyperelliptic curve of genus $g>1$ with the following properties:
\begin{itemize}
\item $X$ has good reduction at 2.
\item $X$ has a rational Weierstrass point $\infty $.
\item At all primes $v$, all points in $X(\Q _v )$ map to a common irreducible component of the special fibre of the minimal regular model of at $v$.
\end{itemize}

Then, as in \cite{QC2}, one constructs a function
\[
h_2 :X(\Q _2 )\to W:=H^1 _f (\Q _2 ,\wedge ^2 V/\Q _2 (1))/\loc _2 H^1 _f (\Q  ,\wedge ^2 V/\Q _2 (1))
\]
with the following properties:
\begin{enumerate}
\item $h_2$ is a nonzero Coleman function. In particular, when restricted to each residue disc it is given by a rigid analytic function (i.e. a convergent power series).
\item $h_2 $ is a quadratic function on $X(\Q )$, in the sense that if $x_i \in X(\Q )$ and $\lambda _i \in \Q _2 $ satisfy 
\[
\sum \lambda _i ([x_i ]-[\infty ])^2 =0
\]
in $\Sym ^2 J(\Q )\otimes \Q $ then 
\[
\sum \lambda _i h_2 (x_i )=0.
\]
Here $J$ denotes the Jacobian of $X$.
\end{enumerate}
Now suppose that the rank of the $J$ is $r$, $X(\Q _2 )_1$ is infinite, and one has found points $x_1 ,\ldots ,x_m \in X(\Q )$, where $m=\frac{r(r+1)}{2}$, such that the dimension of $(x_1 -\infty , \langle h_2 (x_1 )) ,\ldots ,(x_m -\infty , h_2 (x_m ))\rangle $ span $\Sym ^2 J(\Q )\otimes \Q _2$ under the map
\[
P\mapsto P^2 .
\]
Let $\Psi $ denote the map
\[
X(\Q _2 )\to \Sym ^2 H^0 (X_{\Q _2 },\Omega )^* \oplus W.
\]
where the first component sends $z$ to the functional
\[
\omega _0 \omega _1 \mapsto (\int ^z _{\infty }\omega _0 )(\int ^z _{\infty }\omega _1 ),
\]
and the second component is $W$.
Let $W_0 \subset W$ denote the vector space spanned by $\Psi (x_1 ),\ldots ,\Psi (x_m )$. Then $X(\Q _2 )_2$ is equal to the preimage of $W_0$ in $X(\Q _2 )$.

\begin{remark}
In practice, one constructs a generalised height by first constructin a local height (or pre-height)
\[
\widetilde{h}_2 :X(\Q _2 )\to H^1 _f (\Q _2 ,\wedge ^2 V/\Q _2 (1))
\]
with the property that $\widetilde{h}_2 $ is quadratic on $X(\Q )$ (in the sense above) modulo $\loc _2 H^1 _f (\Q ,\wedge ^2 V/\Q _2 (1))$. Hence, to compute $X(\Q _2 )_2$ in the above set-up, one needs $\binom{r(r+1)}{2}+\dim \loc _2 H^1 _f (\Q _2 ,\wedge ^2 V/\Q _2 (1))$ points in $X(\Q )$ which are suitably independent.
\end{remark}
\section{Examples}
Code for the calculations described in this section can be found at \url{https://github.com/netandog/BKdescent}.
\subsection{The BMSST curve}
We consider the example of the curve 
\[
X:y^2 -y =x^5 -x
\]
discussed in the introduction. The $2$-torsion has Galois group $S_5 $. The class groups of the extensions $K:=\Q _f$ and $L:=\Q _f ^{(2)}$ (where $f:=2x^5-32x+16$) are both trivial. 
$X$ has good reduction outside the primes 139 and 449, where it has semistable reduction. We deduce that $H^1 _f (\Q ,\wedge ^2 J[2])$ may be identified with a subspace of 
\[
U:=\Ker (\mathcal{O}_L [1/2]^\times \otimes \F _2 \stackrel{\Nm }{\longrightarrow } \mathcal{O}_K [1/2]^\times \otimes \F _2 ),
\]
which has dimension $5$. Recall that this implies $\dim H^1 _f (\Q ,\wedge ^2 J[2])\leq 4$. Since the Mordell--Weil rank of $J$ is $3$, to deduce finiteness of $X(\Q _2 )_2 $ we need to verify that $\dim H^1 _f (\Q ,\wedge ^2 V)\leq 2$ (note that even if the Mordell--Weil rank was $2$, we would need this to compute $X(\Q _2 )_2$). For any curve $X$ of genus at least one, 
\[
\dim H^1 (\Q ,\wedge ^2 T_2 J)[2]=\dim H^0 (\Q ,\wedge ^2 J[2]) \geq 1
\] 
via the inclusion of Galois modules $\F _2 \hookrightarrow \wedge ^2 J[2]$ dual to the Weil pairing. Hence to prove finiteness of $X(\Q _2 )_2 $ we must prove $\dim H^1 _f (\Q ,\wedge ^2 J[2])=3$.

More precisely, let $\mathfrak{p}$ and $\mathfrak{q}$ denote the unique primes above $2$ in $K$ and $L$ respectively. Then finiteness of $X(\Q _2 )_2$ will follow from finding the subspace 
\[
\mathcal{L}:=H^1 _{f}(\Q _2 ,\wedge ^2 J[2])\subset \Ker (L_{\mathfrak{q} }^\times \otimes \F _2 \stackrel{\Nm }{\longrightarrow } K_{\mathfrak{p}}^\times \otimes \F _2 ),
\]
and verifying
\begin{equation}\label{eqn:verify}
\dim \mathcal{L} \cap \mathrm{Im} (U)\leq 3.
\end{equation}
where $\mathrm{Im} (U)$ is the image of $U$ in $U' :=\Ker (L_{\mathfrak{q} }^\times \otimes \F _2 \stackrel{\Nm }{\longrightarrow } K_{\mathfrak{p}}^\times \otimes \F _2 )$. To do this, we need to compute the crystalline subspace. We do this using our results on integral crystalline lifts and nonabelian $(x-T)$ maps. More precisely, we first computationally verify that $X$ satisfies the conditions of Lemma \ref{lemma:cris_condition}. Hence, if $b$ and $z$ are $\Z _2 $-points of $X-\{ \infty \} $ such that $z-b\in 4\cdot J(\Q _2 )$, then Proposition \ref{prop:main2} defines an element of $U'$ which is a lift of the class $j_{\Gamma ,b}(z)$, and Lemma \ref{lemma:cris_condition} implies that this lies in $\mathcal{L}$ if the image of $X(\Q _2 )$ in $\Sym ^2 H^1 (\Q _2 ,J[2])$ spans $\Sym ^2 H^1 _f (\Q _2 ,J[2])$. This reduces verification of \eqref{eqn:verify} to calculating a basis of $\mathcal{O}_L ^\times $, and calculating whether an element of $L_{\mathfrak{q}}$ is a square. Both of these calculations have been studied in previous work on implementing the Chabauty--Coleman method for hyperelliptic curves (see \cite{stoll:implementing}). We carried out this calculation in magma, searching over $\Z _2 $ points of $X-\{ \infty \} $ to determine a basis of $\mathcal{L}$ and verify \eqref{eqn:verify}.

Having determined the dimension of $\Sel (U_2 )$, we now describe how to compute $X(\Q _2 )_2$. We have the set $S$ of known rational points
\begin{equation}\label{eqn:known}
S=\{ (0,\pm 4),(-2,\pm 4), (2,\pm 4),(\frac{1}{2},\pm \frac{1}{4}),(4, \pm 44), (6, \pm 124), (-\frac{15}{8},\pm \frac{697}{128}),(60,\pm 39436),\infty \}
\end{equation}
By a direct calculation on $J(\Q )$ one may verify that the differences of the symmetric squares of points 
\begin{align*}
(0,4)^{2},(-2,4)^{2}, (2,4)^{2},(\frac{1}{2},\pm \frac{1}{4})^{2},(4, 44)^{2}, \\ 
(6, 124)^{2}, (-\frac{15}{8},\frac{697}{128})^{2},(60,39436)^{2}
\end{align*} 
span $\Sym ^2 J(\Q )\otimes \Q $ (as usual by $P^2$ we mean the class of $(P-\{ \infty \} )\cdot (P-\{ \infty \} )$ in $\Sym ^2 J(\Q )\otimes \Q $). 

For $v=139$ and $449$, a regular semistable model has special fibre an elliptic curve intersecting itself at one point. Hence all $\Q _v$-points lie on a common irreducible component, and we are in the setting above.

Details of the computation of $X(\Q _2 )_2$ may be found in \cite{BKC}. Since the codimension of $\Sel (U_2 )$ in $H^1 _f (\Q _2 ,U_2 )$ is one, in general one would not expect that $X(\Q _2 )_2$ is equal to $X(\Q )$, and indeed we found that $X(\Q _2 )$ contains extra points. Specifically, 
\[
x(X(\Q _2 )_2 ) \subset \left\{  \begin{array}{c}\infty ,0,2,30,1,-1,3,\frac{1}{4},-\frac{15}{16}, \\  24044 +O(2^{15}),22285+O(2^{15}),31171+O(2^{15}),12111+O(2^{15}) \end{array}\right\}.
\]
As in classical implementations of the Chabauty--Coleman method (or in `classical' quadratic Chabauty) one would expect to be able to get around this if one is able to implement the Mordell--Weil sieve (see \cite{BS10}). Roughly speaking, this involves constraining the image of $C(\Q )$ in $J(\Q )/NJ(\Q )$ for a suitably large $N$. Fortunately, for the BMSST curve, Bugeaud, Mignotte, Siksek, Stoll and Tengely have already done this, with \[
N=4449329780614748206472972686179940652515754483274306796568214048000
\]
 (see \cite[\S 10]{BMSST}). They show in particular that any point of $X(\Q )$ must be congruent to the known rational points in the statement of Theorem \ref{thm:BMSST_rational} modulo $2^8$ (with respect to a smooth model over $\Z _2 $). Fortunately, none of the `fake rational points' are congruent to known rational points modulo $2^7$. 
%
%

\begin{remark}
Note that since the size of $S$ is nine, we `only just' have enough rational points to solve for $X(\Q _2 )$, meaning that there are no nontrivial bilinear relations between the generalised heights of rational points. The reader wishing to observe the bilinear properties of the generalised height in practice (or cautious about the accuracy of the computations!) may be interested (or reassured) to know that one can also use points over $\Q (\sqrt{-3})$. Over $\Q (\sqrt{-3})$ we have five additional points
\begin{align*}
& \{(\frac{10+2\sqrt{-3}}{3},\pm \frac{392\sqrt{-3} + 564}{27}), (\frac{-2\sqrt{-3} + 4}{3} ,\pm \frac{-128\sqrt{-3} - 84}{27}), \\
 & (-\sqrt{-3} - 1,\pm (4\sqrt{-3} + 8)),(-\sqrt{-3} + 1,\pm (4\sqrt{-3} + 8)),  (\sqrt{-3} - 3,\pm (4\sqrt{-3} + 32))\}
\end{align*}
(and their Galois conjugates and Weierstrass involutions). Using these points, and known rational points, we can find enough divisors to observe the quadratic properties of $h_2 $.
\end{remark}

\bibliography{bib_BK}

\newcommand{\etalchar}[1]{$^{#1}$}
\begin{thebibliography}{BMS{\etalchar{+}}08}

\bibitem[Bak69]{baker}
A.~Baker.
\newblock Bounds for the solutions of the hyperelliptic equation.
\newblock {\em Proc. Cambridge Philos. Soc.}, 65:439--444, 1969.

\bibitem[BD18]{QC1}
J.~S. Balakrishnan and N.~Dogra.
\newblock Quadratic {C}habauty and rational points, {I}: {$p$}-adic heights.
\newblock {\em Duke Math. J.}, 167, 2018.
\newblock With an appendix by J. Steffen M\"{u}ller.

\bibitem[BD21]{QC2}
Jennifer~S. Balakrishnan and Netan Dogra.
\newblock Quadratic {C}habauty and rational points {II}: {G}eneralised height functions on {S}elmer varieties.
\newblock {\em Int. Math. Res. Not. IMRN}, (15):11923--12008, 2021.

\bibitem[Ber25]{berry}
Lee Berry.
\newblock Refined effective bounds for {B}loch-{K}ato {S}elmer groups associated to hyperelliptic curves.
\newblock {\em arXiv preprint arXiv:2502.11154}, 2025.

\bibitem[BK90]{BK}
S.~Bloch and K.~Kato.
\newblock {L}-functions and {T}amagawa numbers of motives, in {T}he {G}rothendieck {F}estschrift, {V}ol {I}.
\newblock pages 333--400. Birkh\"auser Boston, 1990.

\bibitem[BLLS23]{BLLS}
Dean Bisogno, Wanlin Li, Daniel Litt, and Padmavathi Srinivasan.
\newblock Group-theoretic {J}ohnson classes and non-hyperelliptic curves with torsion {C}eresa class.
\newblock {\em \'{E}pijournal G\'{e}om. Alg\'{e}brique}, 7:Art. 8, 19, 2023.

\bibitem[BLR90]{BLR}
S.~Bosch, W.~Lutkebohmer, and M.~Raynaud.
\newblock N{\'e}ron {M}odels, 1990.

\bibitem[BMS{\etalchar{+}}08]{BMSST}
Yann Bugeaud, Maurice Mignotte, Samir Siksek, Michael Stoll, and Szabolcs Tengely.
\newblock Integral points on hyperelliptic curves.
\newblock {\em Algebra Number Theory}, 2(8):859--885, 2008.

\bibitem[BPS16]{BPS}
Nils Bruin, Bjorn Poonen, and Michael Stoll.
\newblock Generalized explicit descent and its application to curves of genus 3.
\newblock {\em Forum Math. Sigma}, 4:Paper No. e6, 80, 2016.

\bibitem[BS10]{BS10}
Nils Bruin and Michael Stoll.
\newblock The {M}ordell-{W}eil sieve: proving non-existence of rational points on curves.
\newblock {\em LMS J. Comput. Math.}, 13:272--306, 2010.

\bibitem[Cas83]{cassels}
J.~W.~S. Cassels.
\newblock The {M}ordell-{W}eil group of curves of genus {$2$}.
\newblock In {\em Arithmetic and geometry, {V}ol. {I}}, volume~35 of {\em Progr. Math.}, pages 27--60. Birkh\"{a}user, Boston, MA, 1983.

\bibitem[Dog23]{BKdescent2}
Netan Dogra.
\newblock $2$-descent for {B}loch--{K}ato {S}elmer groups {II}.
\newblock {\em preprint}, 2023.

\bibitem[Dog26]{BKC}
Netan Dogra.
\newblock Generalised heights for the {BMSST} curve.
\newblock {\em \url{https://sites.google.com/site/netandogra/links}}, 2026.

\bibitem[FPR94]{FPR94}
J.-M. Fontaine and B.~Perrin-Riou.
\newblock Autour des conjectures de {B}loch et {K}ato: cohomologie galoisienne et valeurs de fonctions {$L$}.
\newblock In {\em Motives ({S}eattle, {WA}, 1991)}, volume~55 of {\em Proc. Sympos. Pure Math.}, pages 599--706. Amer. Math. Soc., Providence, RI, 1994.

\bibitem[Gaz18]{gazaki}
Evangelia Gazaki.
\newblock A finer {T}ate duality theorem for local {G}alois symbols.
\newblock {\em J. Algebra}, 509:337--385, 2018.

\bibitem[HM05]{HM}
Richard Hain and Makoto Matsumoto.
\newblock Galois actions on fundamental groups of curves and the cycle $c-c^-$.
\newblock {\em J. Inst. Math. Jussieu}, 2005.

\bibitem[IM15]{IM15}
Adrian Iovita and Adriano Marmora.
\newblock On the continuity of the finite {B}loch-{K}ato cohomology.
\newblock {\em Rend. Semin. Mat. Univ. Padova}, 134:239--271, 2015.

\bibitem[Kim05]{kim:siegel}
M.~Kim.
\newblock The motivic fundamental group of $\mathbf{P}^1 -\{0,1,\infty\}$ and the theorem of {S}iegel.
\newblock {\em Invent. Math.}, 161(3):629--656, 2005.

\bibitem[Kim09]{kim:chabauty}
M.~Kim.
\newblock The unipotent {A}lbanese map and {S}elmer varieties for curves.
\newblock {\em Publ. Res. Inst. Math. Sci.}, 45(1):89--133, 2009.

\bibitem[KT08]{kim2008component}
M.~Kim and A.~Tamagawa.
\newblock The {$l$}-component of the unipotent {A}lbanese map.
\newblock {\em Math. Ann.}, 340(1):223--235, 2008.

\bibitem[Nek93]{nekovar}
Jan Nekov\'{a}\v{r}.
\newblock On {$p$}-adic height pairings.
\newblock In {\em S\'{e}minaire de {T}h\'{e}orie des {N}ombres, {P}aris, 1990--91}, volume 108 of {\em Progr. Math.}, pages 127--202. Birkh\"{a}user Boston, Boston, MA, 1993.

\bibitem[NSW08]{NSW}
J\"{u}rgen Neukirch, Alexander Schmidt, and Kay Wingberg.
\newblock {\em Cohomology of number fields}, volume 323 of {\em Grundlehren der mathematischen Wissenschaften [Fundamental Principles of Mathematical Sciences]}.
\newblock Springer-Verlag, Berlin, second edition, 2008.

\bibitem[PS97]{PS97}
Bjorn Poonen and Edward~F. Schaefer.
\newblock Explicit descent for {J}acobians of cyclic covers of the projective line.
\newblock {\em J. Reine Angew. Math.}, 488:141--188, 1997.

\bibitem[Sch95]{schaefer}
Edward~F. Schaefer.
\newblock {$2$}-descent on the {J}acobians of hyperelliptic curves.
\newblock {\em J. Number Theory}, 51(2):219--232, 1995.

\bibitem[Ser02]{serre-gc}
J.-P. Serre.
\newblock {\em Galois cohomology}.
\newblock Springer Monographs in Mathematics. Springer-Verlag, Berlin, 2002.
\newblock Translated from the French by Patrick Ion and revised by the author.

\bibitem[SGA72]{SGA7}
{\em Groupes de monodromie en g\'{e}om\'{e}trie alg\'{e}brique. {I}}.
\newblock Lecture Notes in Mathematics, Vol. 288. Springer-Verlag, Berlin-New York, 1972.
\newblock S\'{e}minaire de G\'{e}om\'{e}trie Alg\'{e}brique du Bois-Marie 1967--1969 (SGA 7 I), Dirig\'{e} par A. Grothendieck. Avec la collaboration de M. Raynaud et D. S. Rim.

\bibitem[SS04]{SS}
Edward~F. Schaefer and Michael Stoll.
\newblock How to do a {$p$}-descent on an elliptic curve.
\newblock {\em Trans. Amer. Math. Soc.}, 356(3):1209--1231, 2004.

\bibitem[Sti10]{stix2010trading}
J.~Stix.
\newblock Trading degree for dimension in the section conjecture: the non-abelian {S}hapiro lemma.
\newblock {\em Mathematical Journal of Okayama University}, 52(1), 2010.

\bibitem[Sti13a]{stix2013correction}
J.~Stix.
\newblock Correction to: Trading degree for dimension in the section conjecture: The non-abelian shapiro lemma.
\newblock 2013.

\bibitem[Sti13b]{stix:arithmetic}
Jakob Stix.
\newblock {\em Rational points and arithmetic of fundamental groups}, volume 2054 of {\em Lecture Notes in Mathematics}.
\newblock Springer, Heidelberg, 2013.
\newblock Evidence for the section conjecture.

\bibitem[Sto01]{stoll:implementing}
Michael Stoll.
\newblock Implementing 2-descent for {J}acobians of hyperelliptic curves.
\newblock {\em Acta Arith.}, 98(3):245--277, 2001.

\bibitem[Sto07]{stoll:descent}
Michael Stoll.
\newblock Finite descent obstructions and rational points on curves.
\newblock {\em Algebra Number Theory}, 1(4):349--391, 2007.

\end{thebibliography}
\bibliographystyle{alpha}

\end{document}